\documentclass[preprint,12pt]{elsarticle}

\usepackage[pdftex]{hyperref}
\biboptions{sort&compress,numbers,square}

\usepackage{amsfonts}
\usepackage{graphicx}
\usepackage{epsfig}
\usepackage{geometry}
\usepackage{tikz}
\usetikzlibrary{trees,positioning,shapes}
\usepackage{subcaption}
\usepackage{amsmath}
\usepackage{amsthm}

\usepackage{xcolor, colortbl}

\usepackage{algorithm}
\usepackage{algorithmicx, algpseudocode}

\newcommand{\inred}[1]{\textcolor{red}{#1}}

\usepackage{bm}
\usepackage{mathtools}

\usepackage{rotating}

\begin{document}

\begin{frontmatter}

\title{Deep neural networks for smooth approximation of physics with higher order and continuity B-spline base functions}

\author{Kamil Doleg\l{}o$^{(1)}$, Anna Paszy\'nska$^{(2)}$, Maciej Paszy\'{n}ski$^{(1)}$ and Leszek Demkowicz $^{(3)}$}

\address{$^{(1)}$Institute of Computer Science, \\ AGH University of Science and Technology,
Krak\'{o}w, Poland \\
e-mail: maciej.paszynski@agh.edu.pl \\
$^{(2)}$ Faculty of Physics, Astronomy and Applied Computer Science, \\
Jagiellonian University, Krak\'ow, Poland \\
e-mail: anna.paszynska@uj.edu.pl \\
$^{(3)}$ Oden Institute for Computational and Engineering Sciences, \\
The University of Texas at Austin, USA \\
e-mail: leszek@ices.utexas.edu}

\begin{abstract}
This paper deals with the following important research question. Traditionally, the neural network employs non-linear activation functions concatenated with linear operators to approximate a given physical phenomenon. They "fill the space" with the concatenations of the activation functions and linear operators and adjust their coefficients to approximate the physical phenomena. We claim that it is better to "fill the space" with linear combinations of smooth higher-order B-splines base functions as employed by isogeometric analysis and utilize the neural networks to adjust the coefficients of linear combinations.
In other words, the possibilities of using neural networks for approximating the B-spline base functions' coefficients and by approximating the solution directly are evaluated. 
Solving differential equations with neural networks has been proposed by Maziar Raissi et al. in 2017 \cite{PDESolving} by introducing Physics-informed Neural Networks (PINN), which naturally encode underlying physical laws as prior information. Approximation of coefficients using a function as an input leverages the well-known capability of neural networks being universal function approximators \cite{UniversalApproximators}. In essence, in the PINN approach the network approximates the value of the given field at a given point. We present an alternative approach, where the physcial quantity is approximated as a linear combination of smooth B-spline basis functions, and the neural network approximates the coefficients of B-splines.
This research compares results from the DNN approximating the coefficients of the linear combination of B-spline basis functions, with the DNN approximating the solution directly. We show that our approach is cheaper and more accurate when approximating smooth physical fields.
\end{abstract}
	

\end{frontmatter}

\section{Introduction}
Isogeometric analysis (IGA) has been proposed in 2005 by T.J.R.Hughes et al. \cite{IsogeometricAnalysisProposal} as a generalization of the Finite Element Method into higher order and continuity B-spline basis functions. The method employs smooth B-spline basis to approximate scalar or vector fields representing solutions of different physical phenomena.
The implementation aspects of IGA are summarized in \cite{IGAOvwCompImpl}.
There are a couple of open source libraries offering the isogeometric analysis computational framework. One of them is PetIGA.
PetIGA, a framework for high-performance isogeometric analysis has been developed and described in an article by L. Dalcin, et. al.\cite{PetIGA}. The framework is based
on PETSc \footnote{\url{https://www.mcs.anl.gov/petsc/}}, a high-performance library for the scalable solution of partial differential equations. This library uses traditional solvers for ordinary and partial differential equations. 
Another open source package is GeoPDEs\footnote{\url{http://rafavzqz.github.io/geopdes/}}, a package for isogeometric analysis in MATLAB and Octave.
There are also some MATLAB implementations of stabilized finite element method computations available in an e-book by M. Paszy\'nski\cite{IzogeometrycznaMES}.
Approximation of physics-based fields with splines requires selecting appropriate coefficients for the base functions.  In IGA the coefficients are obtained by solving appropriate systems of linear equations.

The usage of neural networks for solving differential equations has been proposed many times. An article by C. Michoski et al.\cite{DiffEqDNN} presents work on solving partial differential equations with deep neural networks, reviews and extends some of them and focuses on irregular solutions. 

Deep neural networks have recently been shown \cite{NODE} to solve ordinary differential equations with an adaptive precision to speed ratio. In the paper by Yulia Rubinova et. al., neural networks are used to parameterize the derivative of the hidden state instead of the traditional approach of specifying a discrete sequence of hidden layers. In other words, an ODE is embedded into a neural network, by using differential equations solvers as the layers and hidden layers are not predetermined beforehand, but rather the number of layers depends on the desired accuracy. This new family of deep neural network models could possibly replace residual networks. ODE-nets are suitable for time-series data. While ordinary neural networks are discrete and have problems with irregular data, ODE-nets are continuous and allow for evaluation at any point, which results in better accuracy for time-series data. They also have faster testing times at the expense of longer training times. 
Deep Neural Network are also used for optimization of the computational procedures of finite element method \cite{BREVIS2021186,HP}.

An article by M. Raissi et al.\cite{PDESolving} introduced Physics-informed Neural Networks, which naturally encode underlying physical laws as prior information and can be used for both continuous and discrete time models. They are constructed with the help of automatic differentiation to differentiate neural networks with respect to their input coordinates and model parameters. 
%The main disadvantage of this approach is that the general form of the problem itself (eg. a differential equation) must be known beforehand, as it is used in constructing and training the network. Nevertheless, 
PINNs is used for some frequently modelled physical phenomena, like heat transport or flow simulation \cite{nPINNs,RaissiPhysicsIDL}. 

Application of Physics-informed Neural Networks in electric power systems has been recently researched in a preprint by George S. Misyris et. al.\cite{PINNForPowerSystems}. The authors state that the usage of PINNs allows to accurately determine results of differential equations up to 87 times faster than conventional methods. The method used requires less initial training data and can result in smaller neural networks while demonstrating high performance. 

Physics-informed neural networks sometimes fail to be trained, however some research on that topic is being conducted. In the preprint by S. Wang, et al. \cite{WhenPinnsFail} a novel gradient descent algorithm is being proposed to improve PINNs. One of the PINN models' disadvantages over regular deep neural networks, is that they can only predict one specific solution of a Partial Differential Equation. This may be possible to overcome by adapting the PINN learning method to accommodate additional PDE solutions. Nevertheless, they could be trained with minimal to no data from the actual solution if the appropriate boundary conditions were given.

Approximation of coefficients using a function as an input leverages the well-known capability of neural networks being universal function approximators\cite{UniversalApproximators}. This work establishes that the standard multilayer feedforward network architectures can approximate virtually any function of interest to any desired degree of accuracy, provided sufficiently many hidden units are available. It does not, however, address the issue of how many units are needed to attain a given accuracy of approximation. For this paper, it means that a possible failure of the approach with approximating field values could be attributed not to a misuse of neural networks, but rather to inadequate learning, inadequate numbers of hidden units, or the presence of a stochastic rather than a deterministic relation between the input and the target.

This paper proposes a new method to reduce the computational cost of the approximation by applying Deep Neural Networks (DNN) to approximate the base functions' coefficients.  We compare this method to the DNN  approximating the solution obtained directly from IGA solver. We also compare the method to Physics Informed Neural Networks (PINN) \cite{PDESolving} approach where we approximate the field values directly.

There are the following open research questions:

\begin{enumerate}
    % \item \textbf{Which approach is quicker to train?}
    % Methods may differ in terms of training time. Physics-informed neural networks encode laws of physics in their structure 
    \item \textbf{What is the approximation error of similar DNN in both approaches?}
    
    The approximation error may differ for the same problem solved with different approaches. One method might tend to be more accurate than the other methods.
    
    \item \textbf{What is the approximation errors difference between the IGA solution and the solutions from DNN based approaches?} 
    
    If applying DNN in IGA problems proves to be faster than the traditional IGA method, it might still be unusable in practice, because of high approximation errors.  
    
    \item \textbf{Which approach is better?}

    The Universal Approximation Theorem\cite{UniversalApproximators} states that neural networks can represent a variety of functions when given appropriate weights, it does not, however, provide a way to construct those networks and weights, only stating that such a construction is possible.
\end{enumerate}

There are the following research objectives:

\begin{enumerate}
    \item \textbf{Investigating the viability of using DNN for approximating coefficients of base B-spline functions used in IGA.}
    
    Neural networks may prove to be usable in IGA. DNN based approaches should be inspected in terms of accuracy, performance and model size.   
    
    \item \textbf{Comparing the approaches in applying DNN in terms of accuracy and training time.}
    
    Direct solution approximation with DNN, B-spline coefficient approximation using DNN, and PINN approximations may differ in accuracy and training time, due to inherent differences in the nature of all approaches. 
    
%    Finally, the results will be compared with the classical approach.
\end{enumerate}

Paper \cite{PINNIGA} presents an application of PINN into solid mechanics examples, and it compares the PINN model learning based on the exact solution, on the finite element method solution and on the isogeometric analysis solution. It concludes that IGA smooth model results in superior convergence of the PINN training.  In our paper, we show that DNN can actually learn the coefficients of the B-spline basis functions, and it results in a faster convergence and more accurate approximation of the solution than either standard PINN or DNN learning IGA solution.

The structure of our paper is the following. We start in Section 2 by introducing the coefficients of the linear combination of B-splines and their applicability to approximate smooth physical fields. Next, Section 3 compares the three methods employed to incorporate the DNN for the solution of PDEs. Section 4 presents the discussion on the conclusions and future work. We also introduce three Appendixes where we derive the three compared methods on a simple one-dimensional example.

\section{Solving differential equations with neural networks and linear combinations of higher-order and continuity B-spline base functions}

\begin{figure}[ht]
    \centering
    \includegraphics[width=0.9\linewidth]{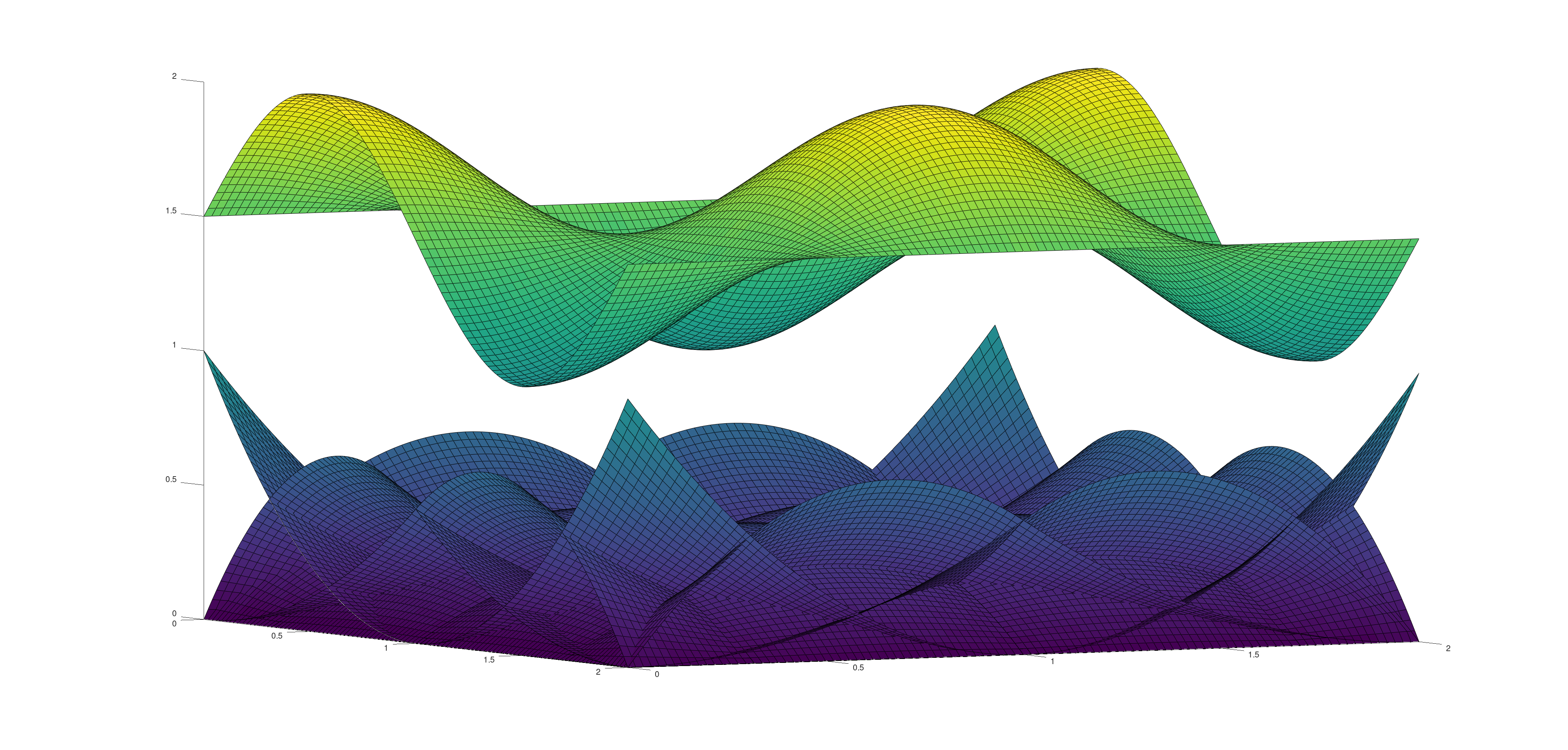}
    % podobnie jak wyżej
    \caption{An example of 2-dimensional quadratic basis generated from $[0,0,0,1,2,2,2]$ knot vectors over $x$ and $y$ axes, and the smooth linear combination of them.}
    \label{fig:2d-base-functions}
\end{figure}

Let assume we want to approximate a function $f(x, y)$ with a linear combination of B-spline functions $u(x, y)$: $$f(x, y) \approx u(x, y)$$
\begin{equation} \label{eq:uxy-2-dim}
u(x,y) = \sum_{i=1, j=1}^{N_x, N_y} u_{i,j} B_i^x B_j^y
\end{equation}
where $B_i^x$, $B_j^y$ denote basis functions over the $x$ and $y$ axis respectively, while $N_x$ and $N_y$ denote numbers of basis functions. An example 2D quadratic base is shown in Figure \ref{fig:2d-base-functions}. The coefficients of the best approximation can be obtained by solving

\begin{flalign}
& \begin{bmatrix}
\int B_1^x B_1^y B_1^x B_1^y dxdy & \hdots & \int B_1^x B_1^y B_{N_x}^x B_{N_y}^ydxdy \\
\vdots & \vdots & \vdots \\
\int B_{N_x}^x B_{N_y}^y B_1^x B_1^y dxdy & \hdots & \int B_{N_x}^x B_{N_y}^y B_{N_x}^x B_{N_y}^y dxdy\\
\end{bmatrix} &
\notag
\\
&
\begin{bmatrix}
u_{1,1} \\
\vdots \\
u_{N_x,N_y}
\end{bmatrix}
= 
\begin{bmatrix}
\int f(x,y) B_1^x B_1^y dxdy \\
\vdots \\
\int f(x,y) B_{N_x}^x B_{N_y}^y dxdy\\
\end{bmatrix}
&
\label{IGA}
\end{flalign}

This system is solvable using traditional solvers but computation-heavy \cite{CostOfContinuity}, despite symmetries in the matrix, especially for a case of non-regular geometry of the computational domain where the B-spline basis functions are span.

The main goal of this project is to verify three different approaches for training deep neural networks in order to provide solutions of IGA simulations with
partial differential equations.
\begin{itemize}
\item{\bf Approximation of coefficients}
The first method is based on "indirect" approximation of the solution, where a DNN predicts the coefficients of the B-spline basis functions that
span over the computational mesh. Linear combinations of those basis functions approximate the solution field.
In this approach, neural networks approximate the values of the base functions' coefficients, $ANN_{i,j}(n)\approx u_{i,j}$, of the linear combination $u_h(x,y) = \sum_{i=1,...,N_x;j=1,...,N_y}u_{i,j} B_i^xB_j^y$ (where $B_i^x$ and $B_j^y$ are the B-spline base functions). This linear combination approximates the exact solution $u_h(x,y)\approx~u(x,y)$ of a given problem. The final solution $u_h$ is obtained by multiplying the relevant base functions by the coefficients predicted by the neural network. There may exist some correlation between the coefficients, that the neural networks can learn, especially in the case of a family of functions.  
The details of this method for simple one-dimensional example are presented in Appendix A.
\item{\bf Direct approximation of the result}
The second method is based on direct approximation of the solution scalar field by a DNN.
Here we have a neural network approximating the solution directly, $ANN(n,x,y) \approx u(x,y)$, we do not use the intermediate linear combination $u_h(x,y)$ in contrast to the previous case.
The details of this method for simple one-dimensional example are presented in Appendix B.
\item{\bf Physics Informed Neural Network}
The third method is based on classical PINN approach \cite{PDESolving}.
The details of this method for simple one-dimensional example are presented in Appendix C.
\end{itemize}

The input data for the DNNs in the first and second cases are the physical model parameters (eg. boundary conditions, material data and forcing), but the outputs from the DNNs black boxes are different in both cases.
In the first case, the output is the coefficients of the linear combination of B-spline basis functions, which depend on the dimensions of the knot vectors and the order and continuity of the basis functions.
In the second case, it is the set of point-wise values of the solution scalar field. The output is thus parameterized with the number of point values along the X and Y-axes.

DNN for the first and second cases learned a family of functions.  Physics-informedneural networks can learn only a selected function from the family, due to the training method used. The DNN for PINN takes as an input the $x$ and $y$ coordinates of the point, and the $n$ parameter is fixed for the PINN neural network.

\section{Comparison of three methods}
\begin{figure}
    \centering
    \includegraphics[width=\linewidth]{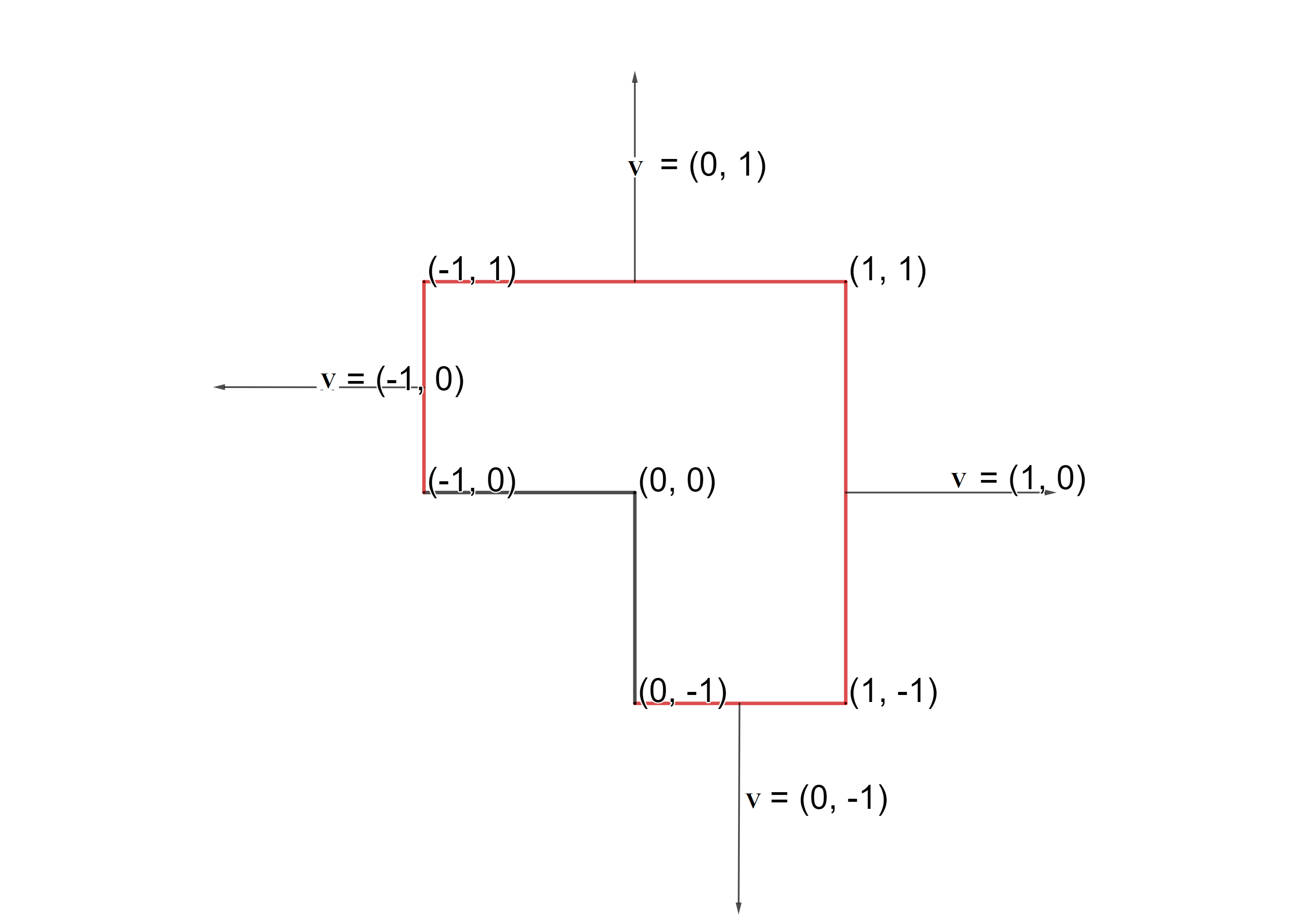}
    \caption{Visualization of the heat transfer problem. The red (external) edge has a Neumann boundary condition $\frac{\partial u}{\partial {\bf v}} (x) = g(x)$ and the black (internal) edge has a Dirichlet boundary condition $u = 0$, where $u$ is the heat function, ${\bf v}$ is the versor normal to the boundary, and $g(x)$ is a given function.}
    \label{fig:heat_problem}
\end{figure}

Heat transfer experiments have been conducted for two-dimensional heat equation on an L-shaped area. External sides are bound with a Neumann boundary condition and internal sides are bound with a Dirichlet boundary condition. This models the heating of the external boundary and fixed zero temperature at the internal boundary.

The equation that we solve is 
$\Delta u = 0$
where $\Delta$ is the Laplacian, namely
$\frac{\partial^2 u}{\partial x^2}+\frac{\partial^2 u}{\partial y^2}=0$
the forcing term is equal to 0, the Dirichlet boundary condition is equal to $0$, and the Neumann boundary condition is 
$\frac{\partial u}{\partial {\bf v}}=g$, 
where we compute the directional derivative in the direction ${\bf v}$ perpendicular to the boundary, and $g(x,y)$ is a given "heating" function.
%We assume the exact solution of this PDE of the form
%$u_{exact}(x,y)=sin(2 \pi n \cdot x) \cdot sin(2 \pi n \cdot y)$.
We enforce the parameterized heating of the external boundary by using a family of "heating" functions, namely
%\begin{eqnarray}
%g(x,y)=2\pi n cos(2 \pi n \cdot x) \cdot sin(2 \pi n \cdot y) \textrm{ for } x=1 \nonumber \\
%g(x,y)=-2\pi n cos(2 \pi n \cdot x) \cdot sin(2 \pi n \cdot y) \textrm{ for } x=-1  \nonumber \\
%g(x,y)=2 \pi n sin(2 \pi n \cdot x) \cdot cos(2 \pi n \cdot y) \textrm{ for } y=1 \nonumber \\
%g(x,y)=-2\pi n sin(2 \pi n \cdot x)  \cdot cos(2 \pi n \cdot y) \textrm{ for } y=-1  \nonumber 
%\end{eqnarray}
\begin{eqnarray}
g(x,y)=v_i2\pi n cos(2 \pi n \cdot x_i) \cdot sin(2 \pi n \cdot x_j) \textrm{ for } |x_i|=1 \nonumber 
\end{eqnarray}
where $i=1,2$, $[ j=(i+1) mod 1]+1$, and ${\bf v}=(v_1,v_2)$ is the versor normal to the boundary.
%We enforce this manufactured solution by computing $g=\frac{\partial u_{exact}}{\partial n}$. In other words, setting g(x,y) in the code for the directional derivative of the manufactured solution makes the code compute the assumed manufactured solution.
By changing the $n$ parameter we adjust the heating function $g$ and we obtain different solutions %over the computational domain 
to the heat transfer problem.

\begin{figure}
    \centering
    \includegraphics[page=1, width=0.4\linewidth]{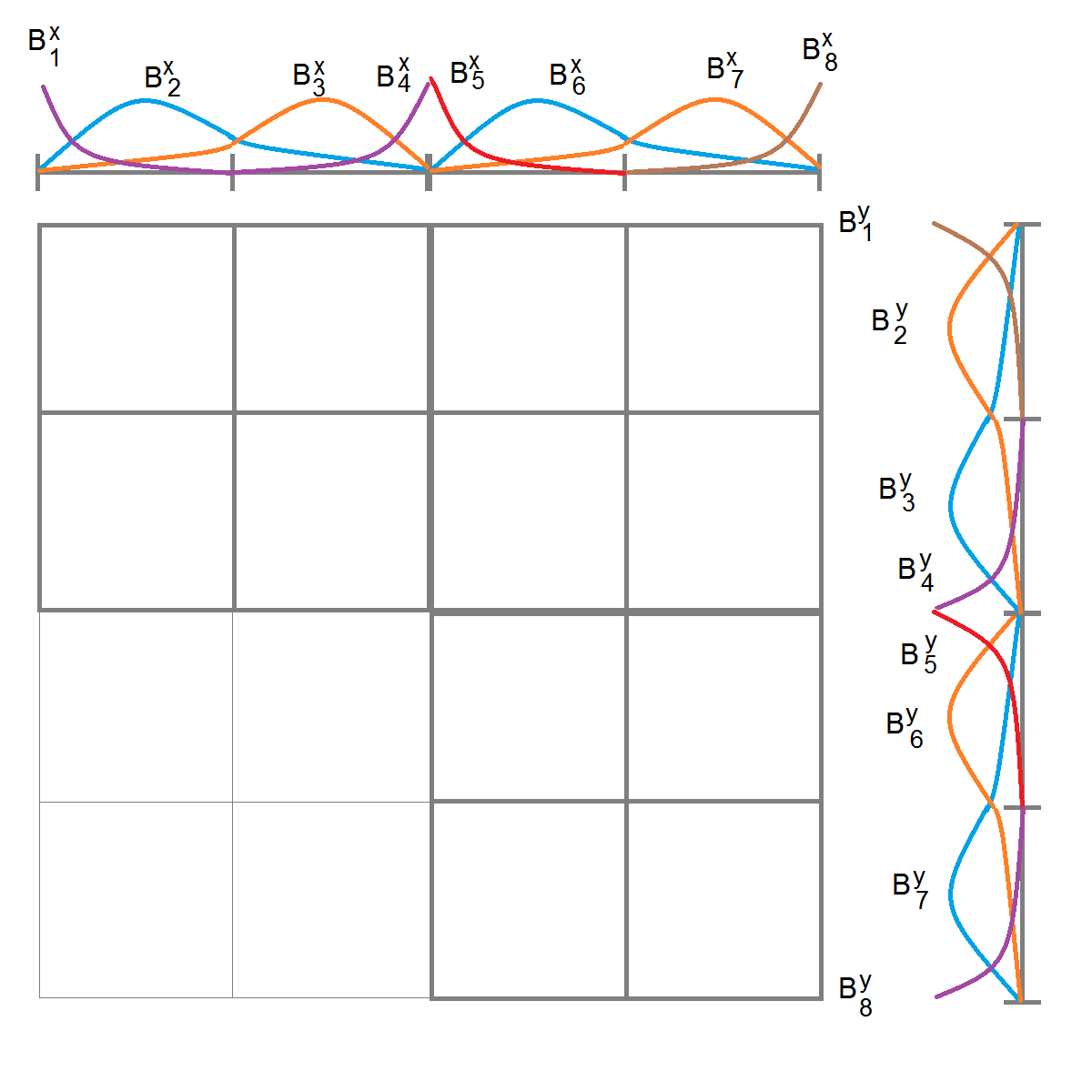}     \includegraphics[page=1, width=0.5\linewidth]{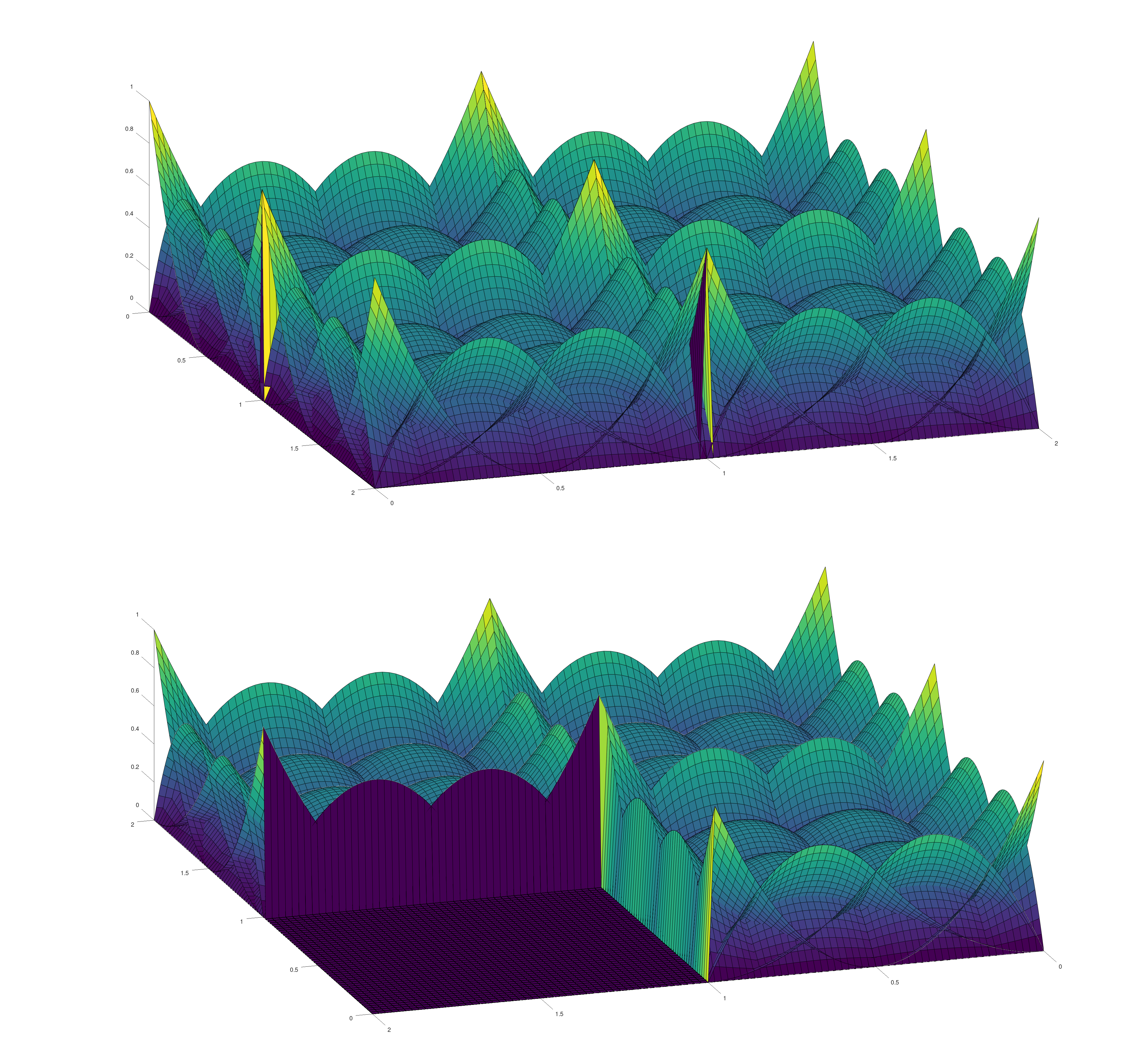} 
    \caption{Knot vectors [0 0 0 0.5 1 1 1 1.5 2 2 2] $\times$ [0 0 0 0.5 1 1 1 1.5 2 2 2] and resulting B-spline base functions, with the splines over left-lower quarter of the domain set to 0.}
    \label{fig:knots}
\end{figure}

Following \cite{IzogeometrycznaMES} (chapter 3), we transform the problem into the weak formulation.
We introduce the domain $\Omega=[-1,1]^2$, 
and we test and integrate by parts to obtain
\begin{eqnarray}
b(u,v)=l(v) \quad \forall v \in V \nonumber \\
b(u,v)=\int_{\Omega} \nabla u \cdot \nabla v dxdy \nonumber \\
l(v)=\int_{\partial \Omega} g v dS
\end{eqnarray}
We discretize with B-spline base functions.
We introduce the knot vectors along $x$ and $y$ axis of coordinates, where we repeat the knot at 0, to obtain the tensor product basis function $B^x_iB^y_j$, $i=1,...,N_x$, $j=1,...,N_y$. 
An example of knot vectors [0 0 0 0.5 1 1 1 1.5 2 2 2] $\times$ [0 0 0 0.5 1 1 1 1.5 2 2 2] and resulting B-spline base functions with B-splines over lower-left corner of the domain set to 0, are illustrated in Figure \ref{fig:knots}. 

The coefficients of B-splines employed for training of the DNN are obtained by solving the system of equations

\begin{flalign}
& \begin{bmatrix}
b\left(B_1^x B_1^y, B_1^x B_1^y\right) & \hdots & b\left( B_1^x B_1^y, B_{N_x}^x B_{N_y}^y \right)\\
\vdots & \vdots & \vdots \\
b\left( B_{N_x}^x B_{N_y}^y, B_1^x B_1^y \right) & \hdots & b\left( B_{N_x}^x B_{N_y}^y, B_{N_x}^x B_{N_y}^y\right) \\
\end{bmatrix} &
\notag
\\
&
\begin{bmatrix}
u_{1,1} \\
\vdots \\
u_{N_x,N_y}
\end{bmatrix}
= 
\begin{bmatrix}
l\left(B_1^x B_1^y \right) \\
\vdots \\
l\left( B_{N_x}^x B_{N_y}^y \right)\\
\end{bmatrix}
&
\label{IGA}
\end{flalign}

Additionally, for the rows related to B-splines located in the left-lower quarter of the domain, we set the rows to zero, we put 1.0 on the diagonal and zero on the right-hand side.
An alternative way of setting boundary conditions would be to add an additional non-trainable layer to the neural network enforcing the Dirichlet boundary condition following the ideas described in \cite{bc}.

\subsection{Coefficient approximation  by DNN} \label{sub:coefficient-approximation}

We'd like the network to learn a family of solutions of the heat transfer problem.
%form $u_{exact}(x,y)=sin(2 \pi n \cdot x) \cdot sin(2 \pi n \cdot y)$. 
The input to the network is the $n$ parameter. 
The output to the network are the coefficients $u_{ij}$ of the linear combination of B-splines.
%The number of coefficients is equal to the number of the network's output neurons, one neuron is responsible for one coefficient value; 

\begin{figure}
    \centering
    \includegraphics[page=1, width=\linewidth]{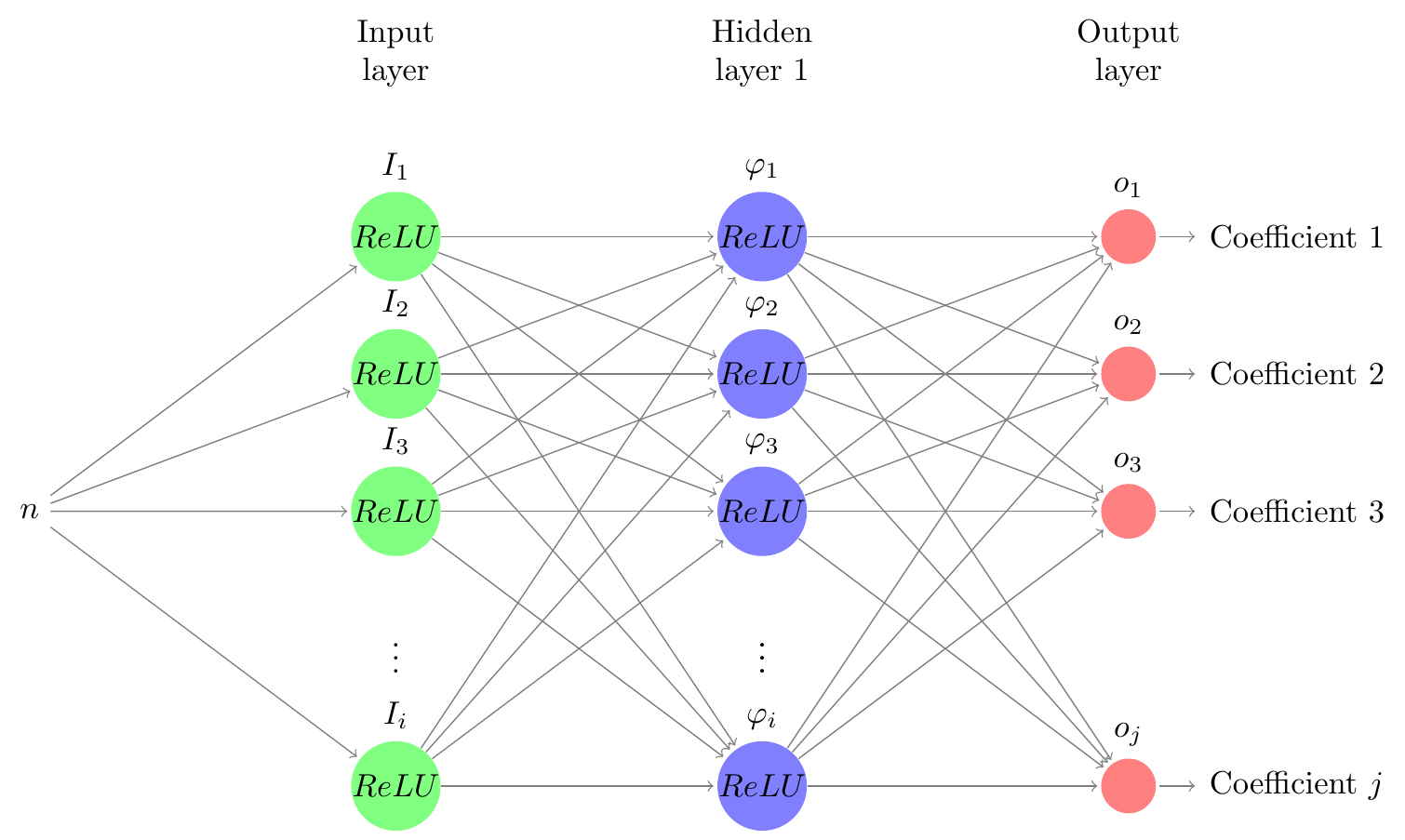} 
    \caption{Visualization of the neural network architecture. The value of $i$ is between 100 and 1000, and the value of $j$ is equal to the number of coefficients, as explained in Table \ref{tab:heat-dnn-coeff-arch}. }
    \label{fig:heat-dnn-coeff-arch-vis}
\end{figure}

The network architecture varies with input data - the output layer has to have the same number of neurons as coefficients. This is the main drawback of this approach, as the network has to be reconfigured with each mesh change. The model was trained with a learning rate reduction on plateau. 
\begin{table}[H]
\centering
\begin{tabular}{|l|l|l|}
\hline
Layer        & Number of neurons & Activation \\ 
& & function \\ \hline
input        & 100 - 1000    & ReLU          \\ \hline
hidden layer 1      & 100 - 1000 & ReLU          \\ \hline
output 2     & equal to & none          \\
 & the number of coefficients   & \\ \hline
\end{tabular}
    \caption{Network architecture of the DNN used to approximate the coefficients of the heat transfer solution. The number of neurons in the hidden layer varied with inputs. For 169 output coefficients the hidden layer had 100 neurons and the number rose to as much as 1000 neurons for 2025 output coefficients. ReLU stands for Rectified Linear Unit.}
    \label{tab:heat-dnn-coeff-arch}
\end{table}

%\subsubsection{Results} 

% 41.88719415664673 seconds for heat_coefficients_01_09_002_quadratic_10el_nx_ny

\begin{figure}[H]
    \centering
    \includegraphics[width=0.49\linewidth]{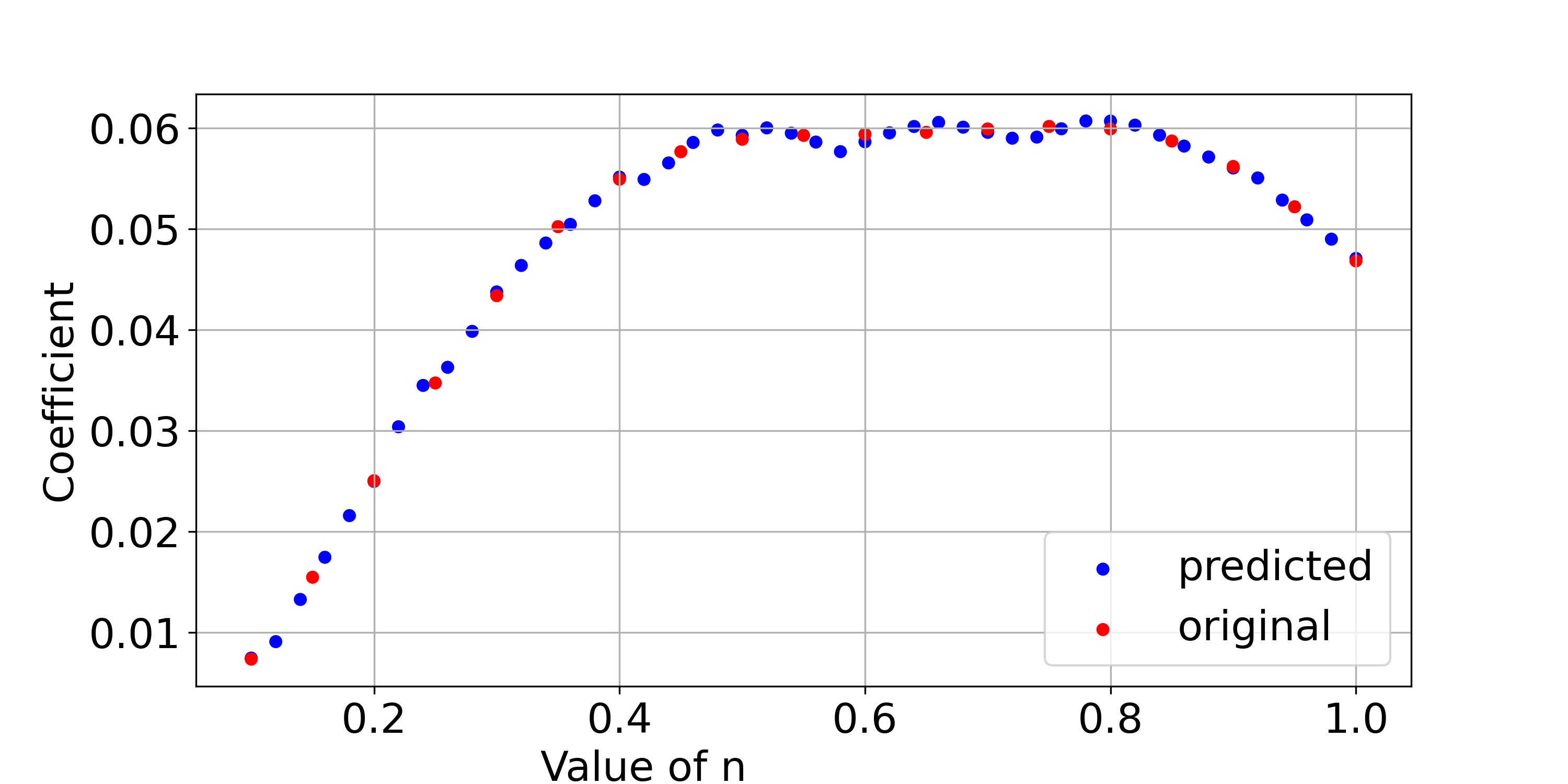} \includegraphics[width=0.49\linewidth]{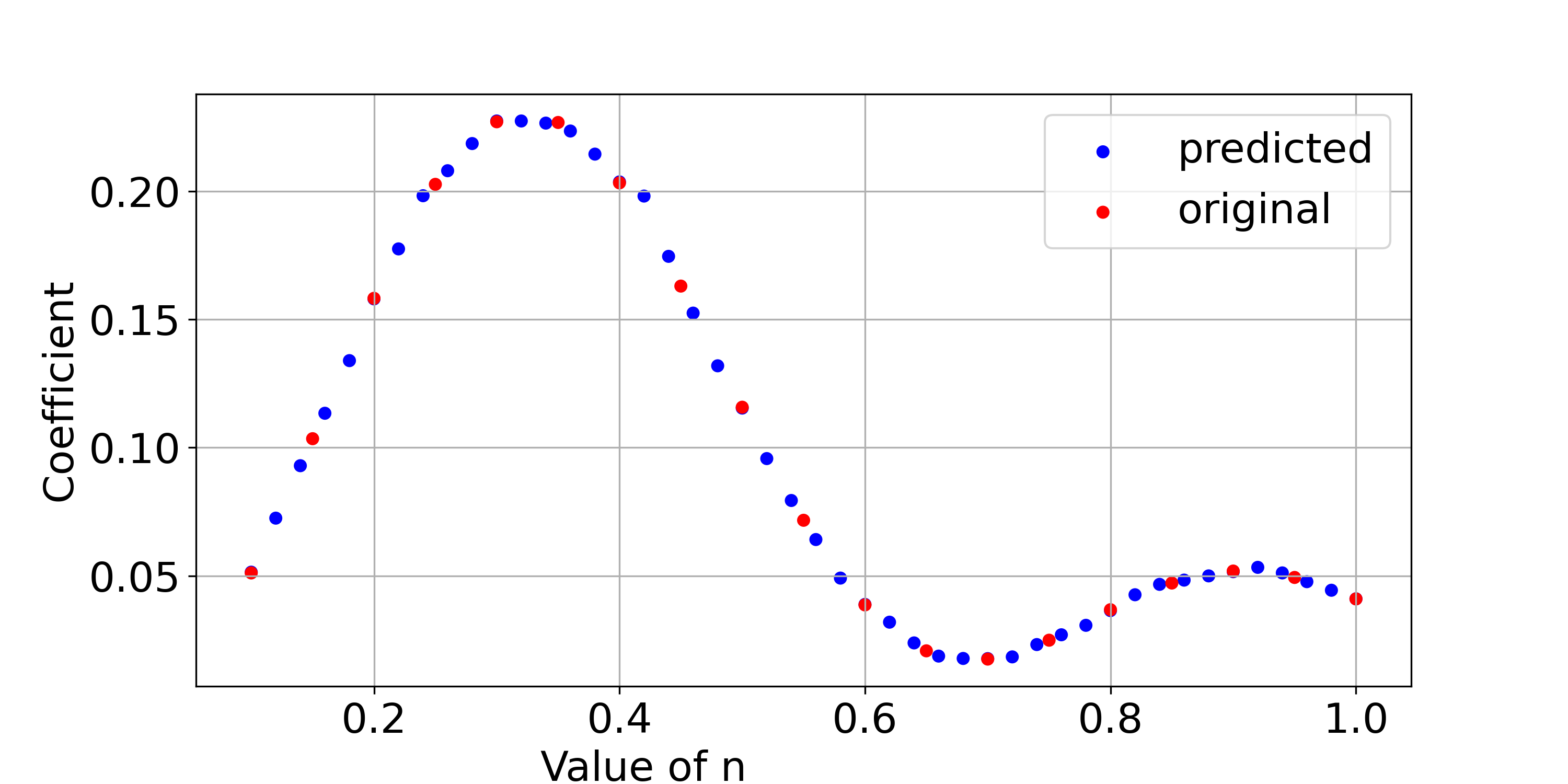}

    \caption{Approximating the $B_{5,2}(x)B_{4,2}(y)$ and $B_{7,2}(x)B_{7,2}(y)$ coefficients with a neural network trained to approximate the coefficient value directly}
    \label{fig:heat-simple_coef6-direct}
\end{figure}

\begin{figure}[H]
    \centering
    \includegraphics[width=0.49\linewidth]{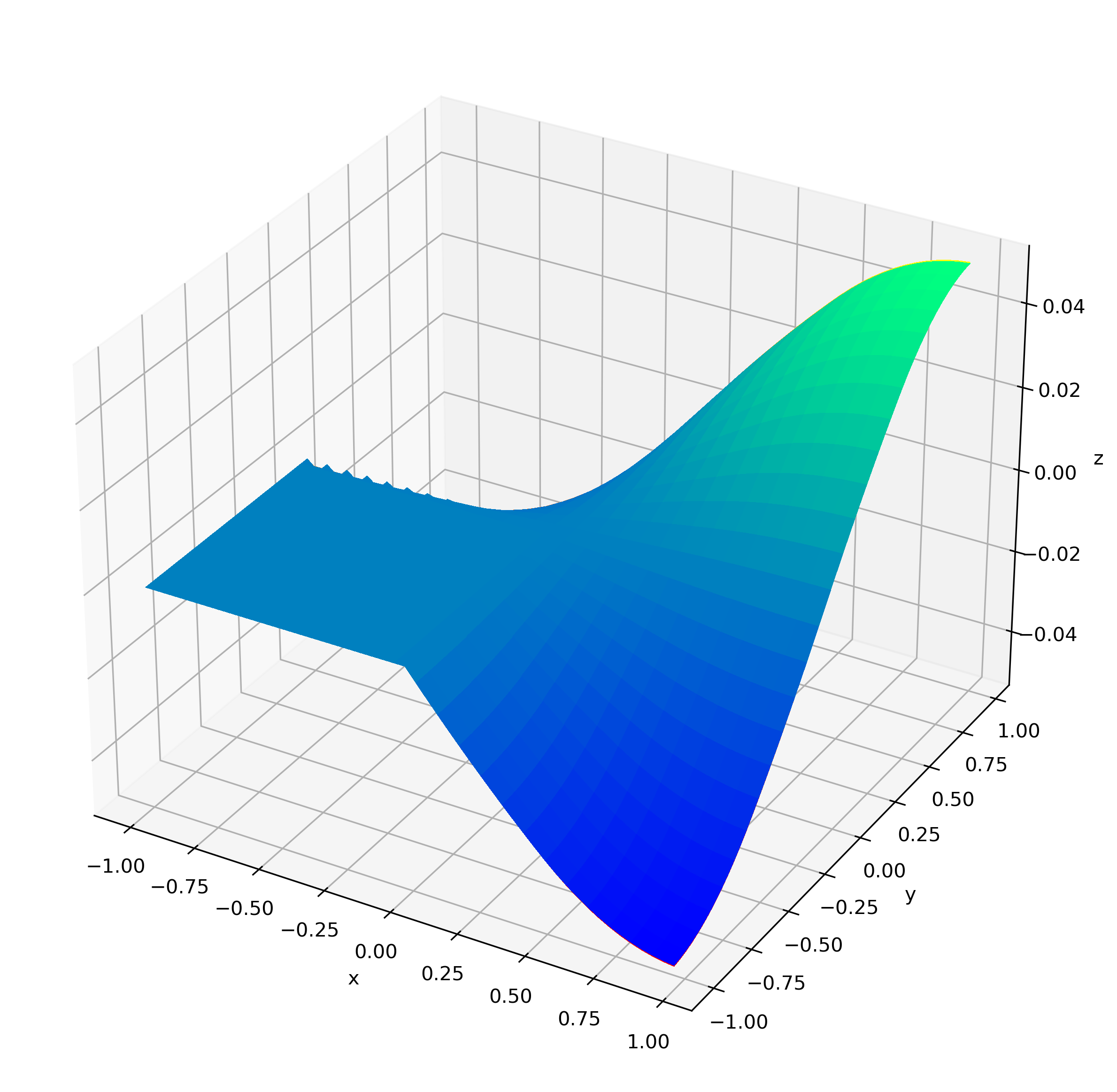}    \includegraphics[width=0.49\linewidth]{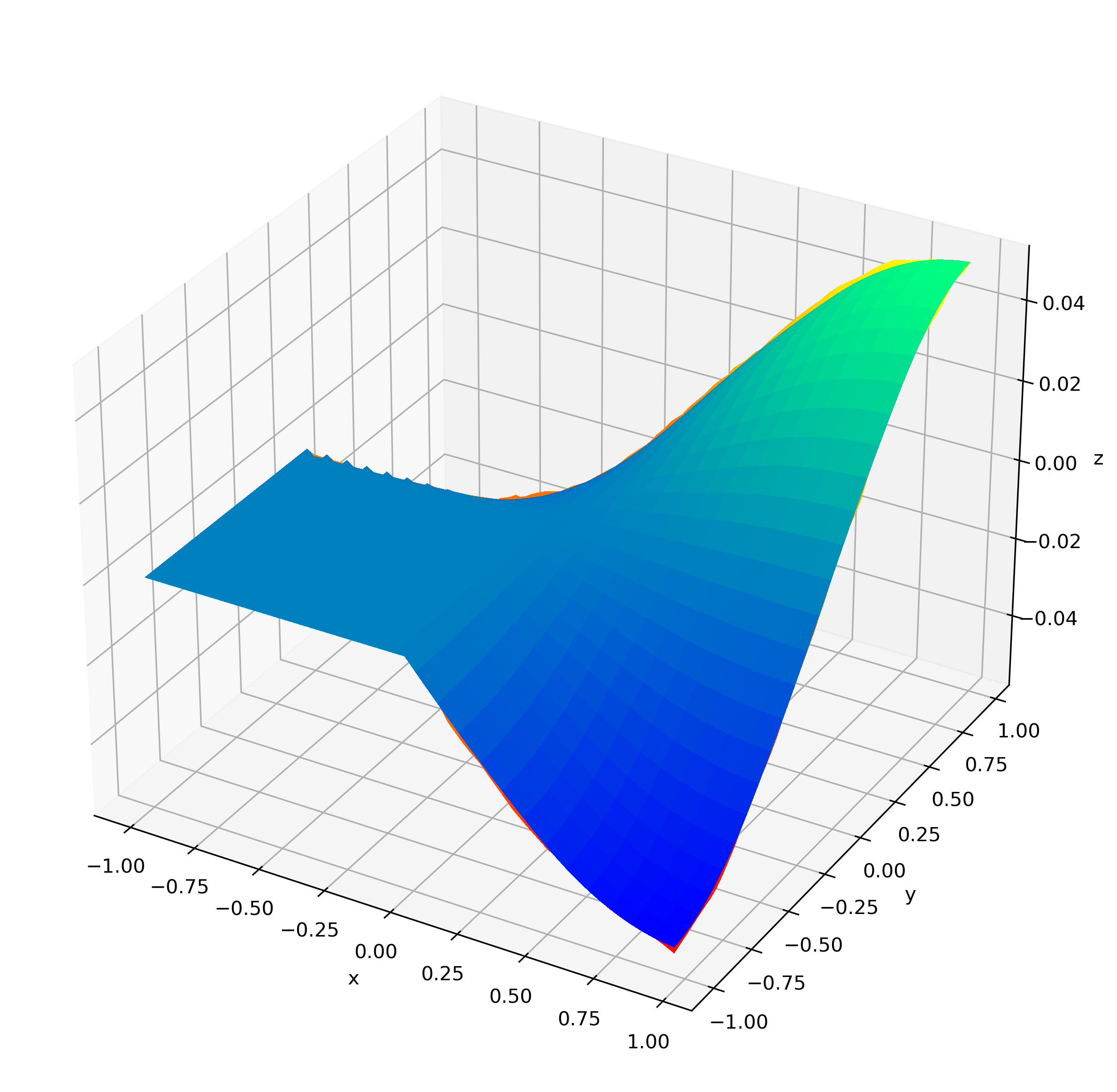}
    \caption{{\bf Left panel:} The IGA solution of (\ref{IGA}) (green-blue ) and the coefficients approximating DNN (orange-red) solutions. The solutions are practically identical. The MSE of the coefficients is $2.27\text{e-}8$, the pointwise MSE is $6.60\text{e-}9$. {\bf Right panel:} The IGA solution of (\ref{IGA}) (blue-green) and the DNN approximating the solution directly (red-yellow). The solutions are practically identical. The mean squared error is approximately $8.42\text{e-}7$.}

    \label{fig:fem-dnn-solution}
\end{figure}

We present results obtained for 7 quadratic B-spline functions in each direction spanned over the $[0, 0, 0, 1, 2, 2, 3, 4, 4, 4]$ knot vector, which gives $7 \times 7 = 49$ coefficients in total.
Coefficient approximation yielded satisfactory results, see left panel in Figure \ref{fig:fem-dnn-solution}. 

\begin{figure} [H]
    \centering
    \includegraphics[width=0.8\linewidth]{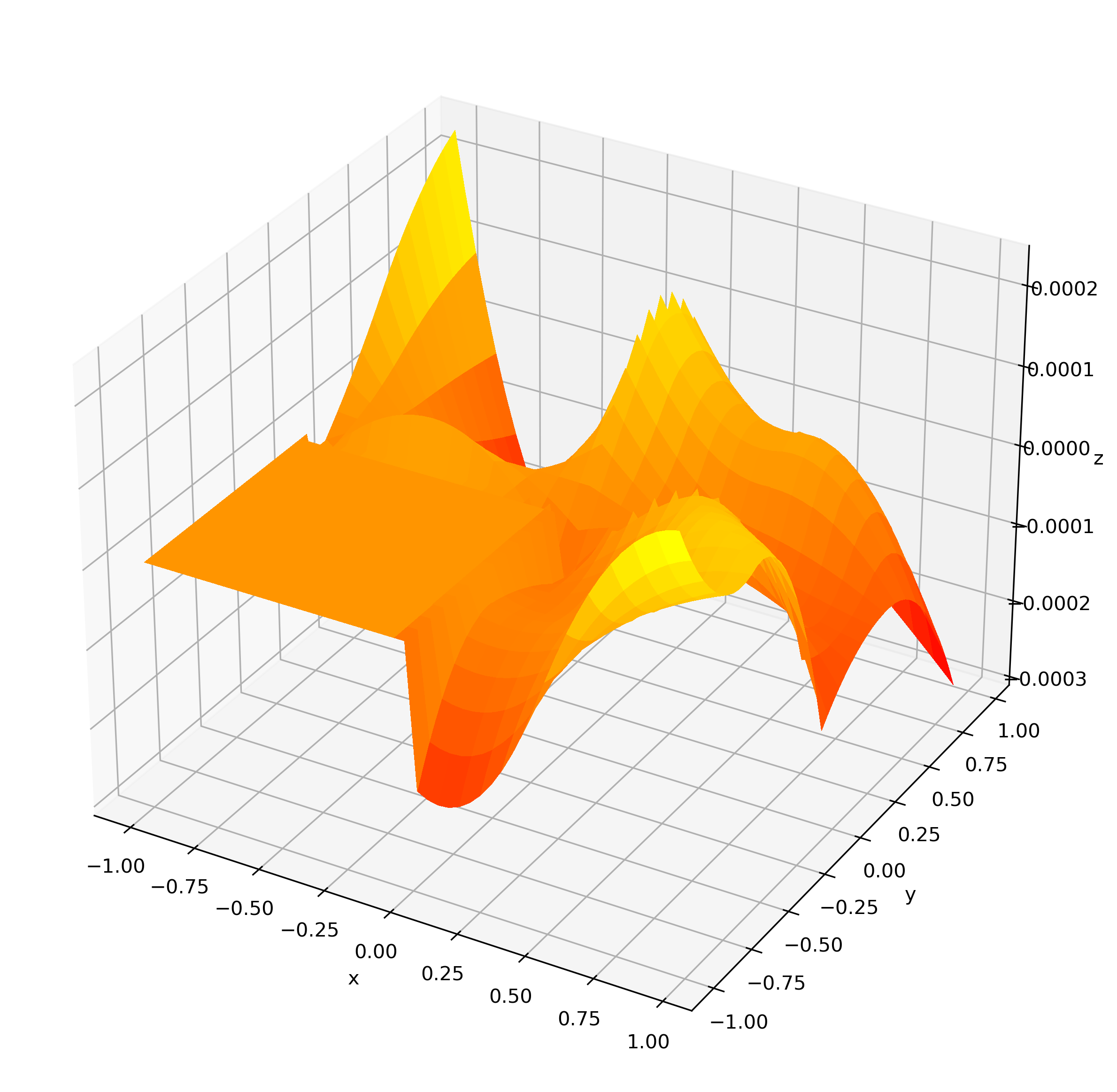}
    \caption{Difference between the solutions constructed using coefficients from IGA solver (\ref{IGA}) and the coefficients obtained from trained DNN. The mean squared error is approximately $6.60\text{e-}9$.}
    \label{fig:error-smooth-PINN}
\end{figure}

\begin{figure}[H]
    \centering
    \includegraphics[width=\linewidth]{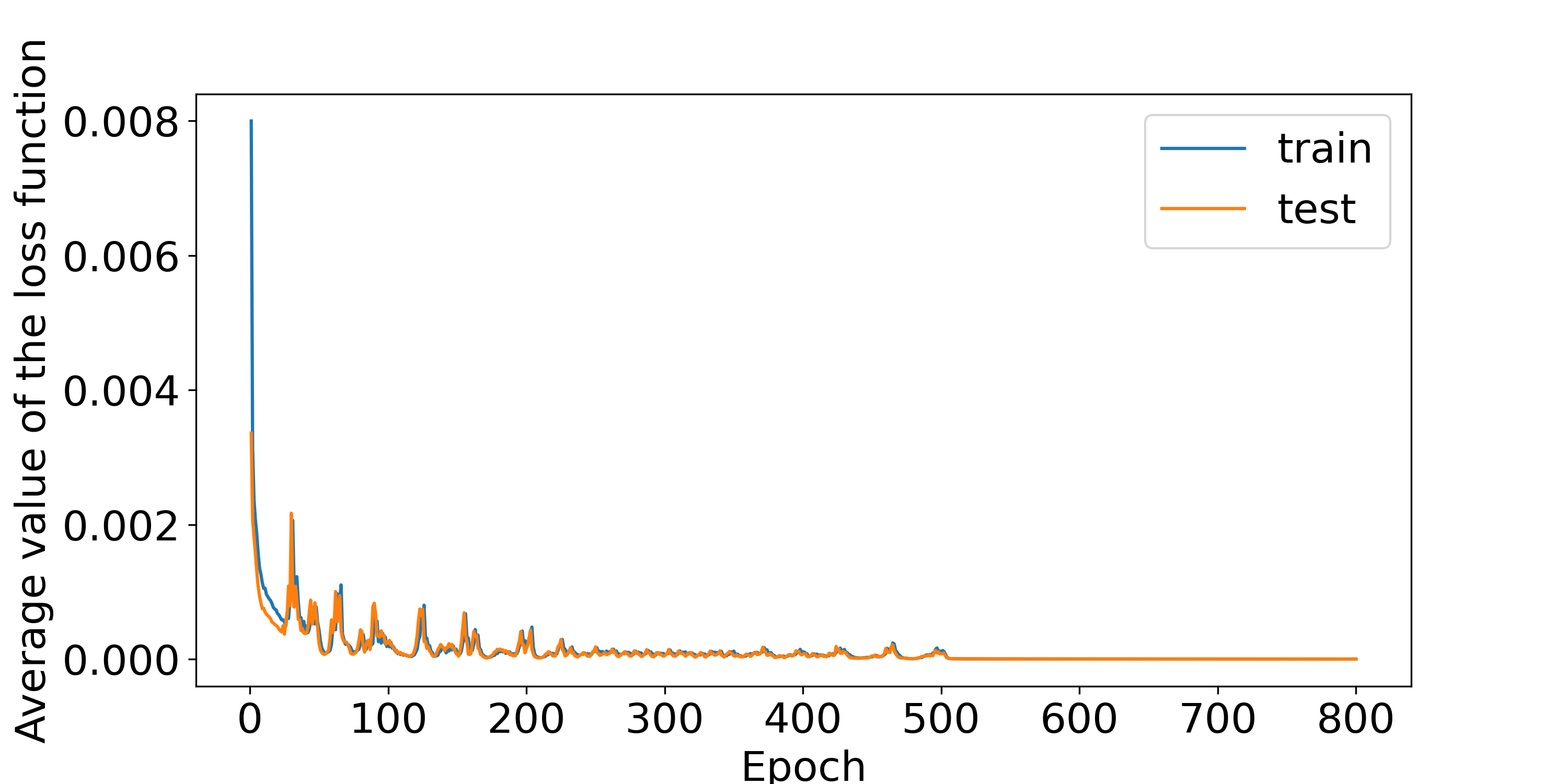}
    \caption{Mean squared error loss function averaged over an epoch for the first 800 out of 2 000 epochs. Each batch contained 19 entries.}
    \label{fig:heat-dnn-coeff-loss}
\end{figure}

The neural network could learn a family of heat functions and accurately predict the coefficients for the base functions, as shown in Figure \ref{fig:heat-simple_coef6-direct}. 
The difference between the solution obtained with IGA solver (\ref{IGA}) and the solution obtained from the DNN trained to approximate the coefficients is shown in the Figure \ref{fig:error-smooth-PINN}. Model convergence is visible in the Figure \ref{fig:heat-dnn-coeff-loss}. Table \ref{tab:heat-dnn-coeff-results} summarizes the results for several different meshes and different base functions.
\begin{table}[H]
\centering
\begin{tabular}{|l|l|l|l|l|}
\hline
Mesh  & \# coeffs & Spline degree & MSE value & Training \\ 
 & & & & time [s] \\ \hline
10 x 10 & 121 & linear & 1.98e-07 & 9  \\ \hline
10 x 10 & 169 & quadratic & 5.29e-07 & 10 \\ \hline
10 x 10 & 225 & cubic & 2.18e-07 & 11  \\ \hline
20 x 20 & 441 & linear & 6.99e-07 & 14  \\ \hline
20 x 20 & 529 & quadratic & 6.96e-07 & 18  \\ \hline
20 x 20 & 625 & cubic & 9.14e-08 & 19 \\ \hline
40 x 40 & 1681 & linear & 1.46e-06 & 15 \\ \hline
40 x 40 & 1849 & quadratic & 1.37e-06 & 16 \\ \hline
40 x 40 & 2025 & cubic & 7.58e-06 & 18  \\ \hline
% 10 x 10 & 121 & linear & 1.98e-07 & 9.85  \\ \hline
% 10 x 10 & 169 & quadratic & 5.29e-07 & 10.28  \\ \hline
% 10 x 10 & 225 & cubic & 2.18e-07 & 11.05  \\ \hline
% 20 x 20 & 441 & linear & 6.99e-07 & 14.67  \\ \hline
% 20 x 20 & 529 & quadratic & 6.96e-07 & 18.24  \\ \hline
% 20 x 20 & 625 & cubic & 9.14e-08 & 19.66 \\ \hline
% 40 x 40 & 1681 & linear & 1.46e-06 & 15.30 \\ \hline
% 40 x 40 & 1849 & quadratic & 1.37e-06 & 16.30 \\ \hline
% 40 x 40 & 2025 & cubic & 7.58e-06 & 18.37  \\ \hline
\end{tabular}
    \caption{Summary of obtained training results for a coefficient-approximating DNN. The training time rises with the number of coefficients, as expected. The training times are only an indication of the time required to train the model, because they vary significantly with batch size and the number of iterations. The error tends to increase with the number of coefficients, however the actual output is still perfectly usable.}
    \label{tab:heat-dnn-coeff-results}
\end{table}

\subsection{DNN approximating the solution directly} \label{sub:direct-approximation-with-dnn}
\begin{figure} [H]
    \centering
    \includegraphics[width=0.8\linewidth]{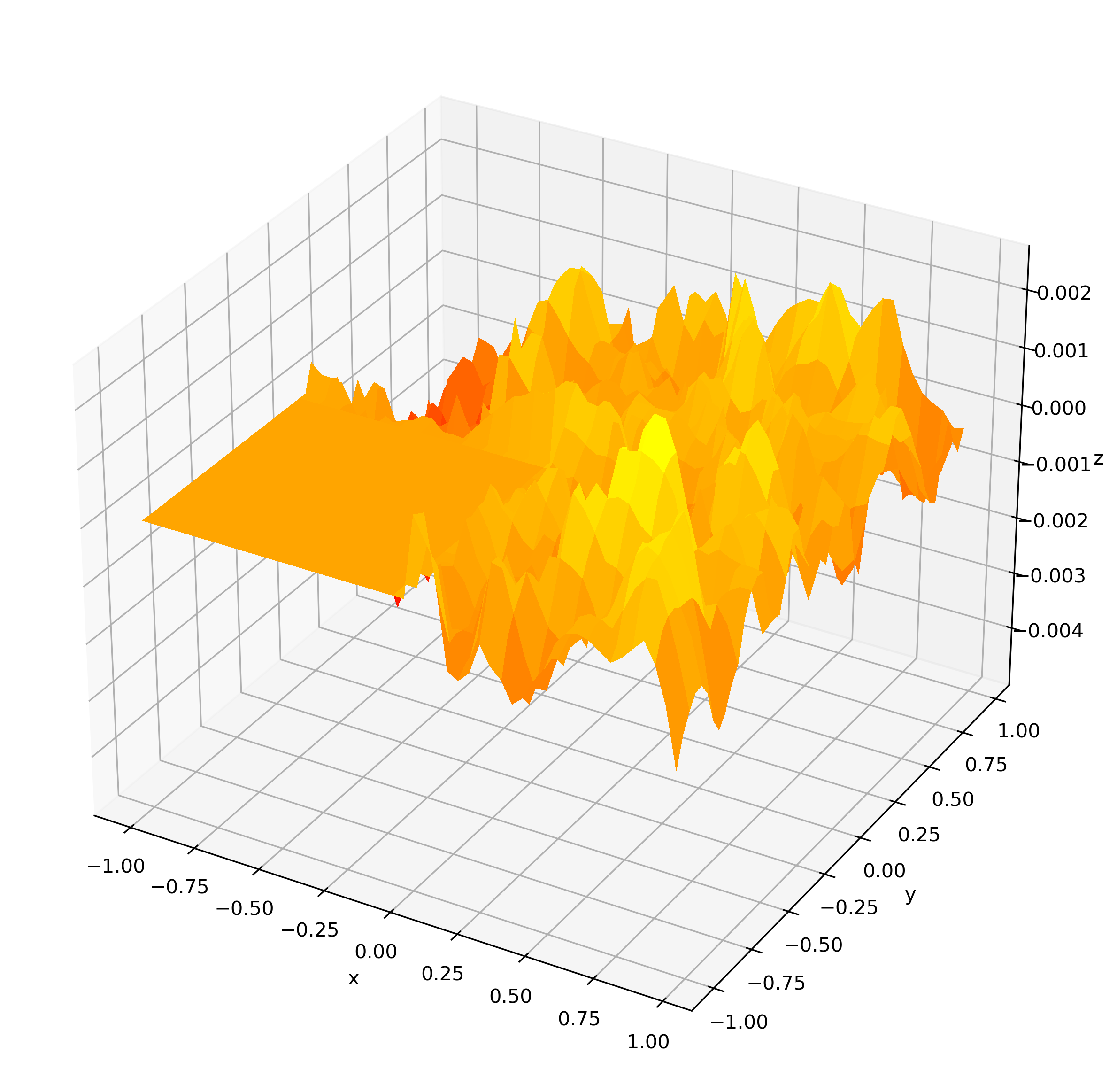}
    \caption{Difference between the solution obtained with IGA solver (\ref{IGA}) and the the DNN approximating the solution directly.
The mean squared error is approximately $8.42\text{e-}7$.}
    \label{fig:error-PINN}
\end{figure}

\begin{figure}[H]
    \centering
    \includegraphics[width=\linewidth]{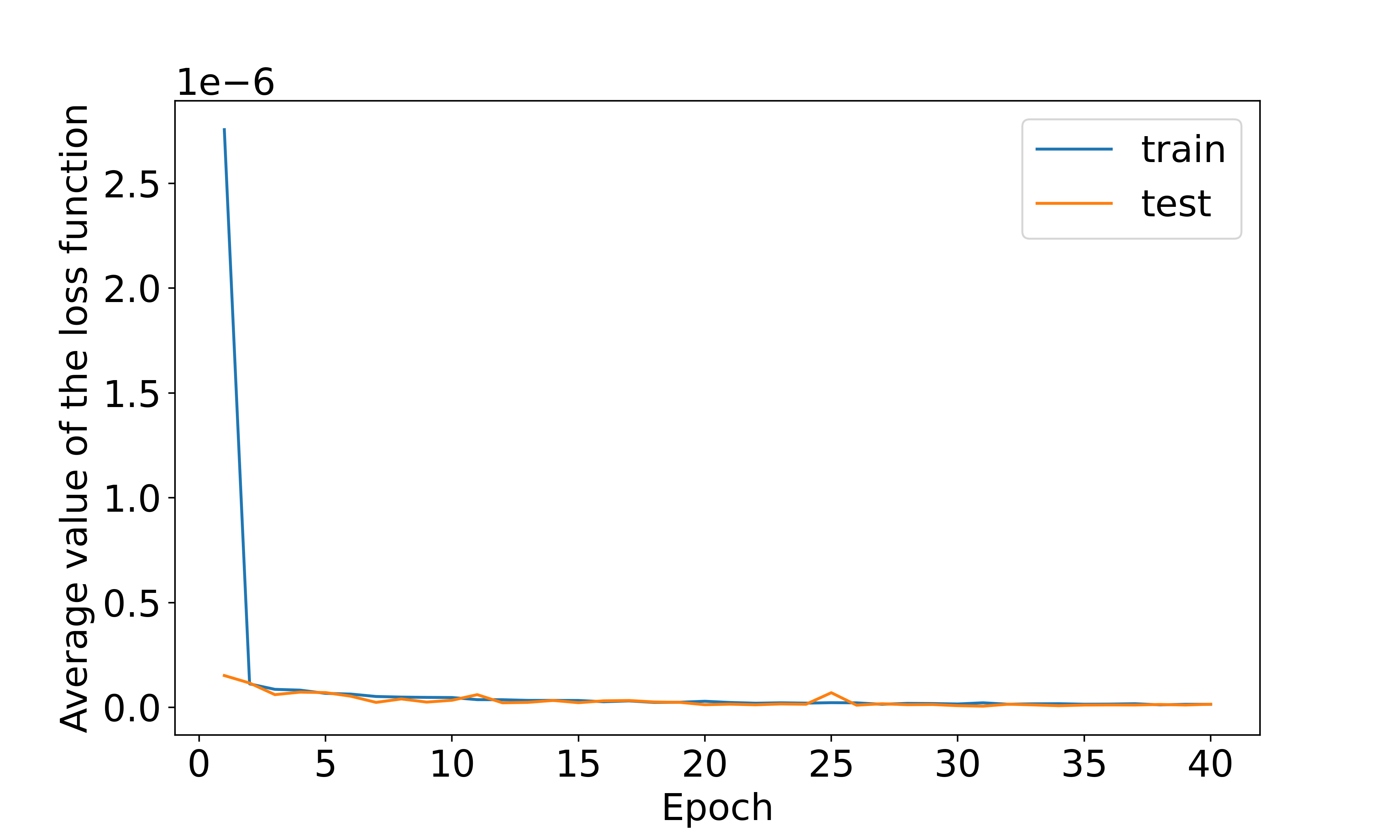}
    \caption{Mean squared error loss function over 40 training epochs.}
    \label{fig:heat-dnn-direct-loss}
\end{figure}

A deep neural network with 2 hidden layers has been trained to directly approximate the result. 
The network was designed to take 3 arguments - the value of $n$, and $x$ and $y$ coordinates of the desired point of the solution. The model was a fully-connected feed-forward neural network with 2 hidden layers, 100 neurons each. 
Direct approximation also yielded satisfactory results. The neural network could learn a family of heat functions and accurately predict the solutions, as shown in right panel in Figure \ref{fig:fem-dnn-solution}. 
The solution is slightly less accurate than the one obtained with coefficient approximation, however it should still be usable. A difference between the solution obtained from IGA solver (\ref{IGA}) and the DNN approximating the solution directly is shown in Figure \ref{fig:error-PINN}.
Model convergence is visible in the Figure \ref{fig:heat-dnn-direct-loss}. 
Table \ref{tab:heat-times-comparison} summarizes the learning times for several different meshes and different base functions and compares them to the learning times of the $n$ input network. It is worth mentioning here that the network approximating the solution does not have to be retrained for different meshes, in contrast to the coefficient approximation, which has to be retrained for different meshes, so the training times here should be treated only as an indication of the order of magnitude.

\begin{table}[H]
\centering
\begin{tabular}{|l|l|l|}
\hline
Layer        & Number of neurons & Activation function \\ \hline
input        & 100    & ReLU          \\ \hline
hidden layer 1      & 100  & ReLU          \\ \hline
hidden layer 2      & 100  & ReLU          \\ \hline
output       & 1    & none          \\ \hline 
\end{tabular}
    \caption{Summary of the layers and activation functions of the DNN used to directly approximate the heat transfer solution. ReLU stands for Rectified Linear Unit.}
    \label{tab:heat-dnn-direct-arch}
\end{table}

\begin{figure}[H]
    \centering
    \includegraphics[page=3, width=\linewidth]{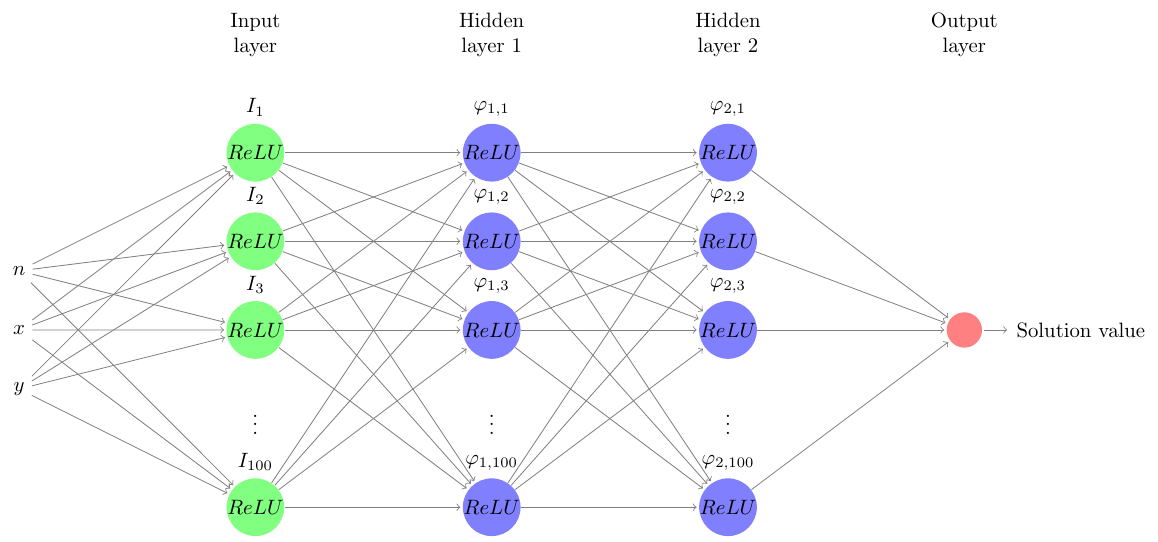} 
    \caption{Visualization of the architecture of the DNN used to directly approximate the heat transfer solution.}
    \label{fig:heat-dnn-direct-arch-vis}
\end{figure}

%\subsubsection{Results} 

\begin{table}[H]
\centering
\begin{tabular}{|l|l|l|l|}
\hline
Mesh size & Spline degree & Coefficient & Direct \\ 
& & network  & network \\ 
          &               & training & training  \\
          &               & time [s]    & time [s] \\ \hline
10 x 10 & linear & 9 & 115  \\ \hline
10 x 10 & quadratic & 10 & 114  \\ \hline
10 x 10 & cubic & 11 & 48  \\ \hline
20 x 20 & linear & 14 & 43  \\ \hline
20 x 20 & quadratic & 18 & 47  \\ \hline
20 x 20 & cubic & 19 & 42  \\ \hline
40 x 40 & linear & 15 & 44  \\ \hline
40 x 40 & quadratic & 16 & 42  \\ \hline
40 x 40 & cubic & 18 & 48  \\ \hline
% 10 x 10 & 9.85 & 115.08  \\ \hline
% 10 x 10 & 10.28 & 114.45  \\ \hline
% 10 x 10 & 11.05 & 48.87  \\ \hline
% 20 x 20 & 14.67 & 43.02  \\ \hline
% 20 x 20 & 18.24 & 47.74  \\ \hline
% 20 x 20 & 19.66 & 42.68  \\ \hline
% 40 x 40 & 15.30 & 44.14  \\ \hline
% 40 x 40 & 16.30 & 42.40  \\ \hline
% 40 x 40 & 18.37 & 48.19  \\ \hline
\end{tabular}
    \caption{Comparison of training times of the coefficient-approximating neural network and the solution approximating neural network.% The training times are only an indication of the time required to train the model, because they vary significantly with batch size and the number of iterations. 
The MSE of achieved outputs had the same order of magnitude.}
    \label{tab:heat-times-comparison}
\end{table}

\subsection{Direct approximation with PINN}
Deep neural networks in previous sections  learned a family of functions. Physics-informed neural networks can learn only a selected function from the family, due to the training method used.
The network was designed to take 2 arguments - the value of $x$ and $y$ coordinates of the desired point of the solution. The model was a fully-connected feed-forward neural network with 2 hidden layers, 50 neurons each: 

\begin{table}[H]
\centering
\begin{tabular}{|l|l|l|}
\hline
Layer        & Number of neurons & Activation function \\ \hline
input        & 50  & ReLU          \\ \hline
hidden layer 1      & 50  & ReLU          \\ \hline
hidden layer 2      & 50  & ReLU          \\ \hline
output       & 1   & none          \\ \hline 
\end{tabular}
    \caption{Summary of the layers and activation functions of the PINN used to directly approximate the heat transfer solution. ReLU stands for Rectified Linear Unit.}
    \label{tab:heat-pinn-direct-arch}
\end{table}

\begin{figure}[H]
    \centering
    \includegraphics[page=4, width=\linewidth]{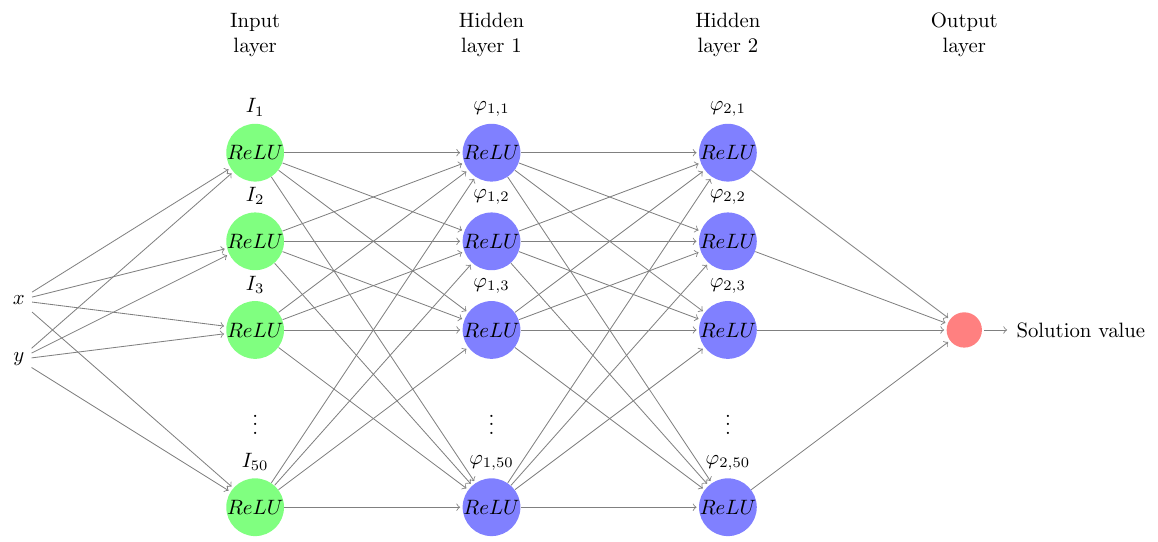} 
    \caption{Visualization of the architecture of the PINN used to directly approximate the heat transfer solution.}
    \label{fig:heat-pinn-direct-arch-vis}
\end{figure}

\begin{figure}[H]
    \centering
    \includegraphics[width=0.8\linewidth]{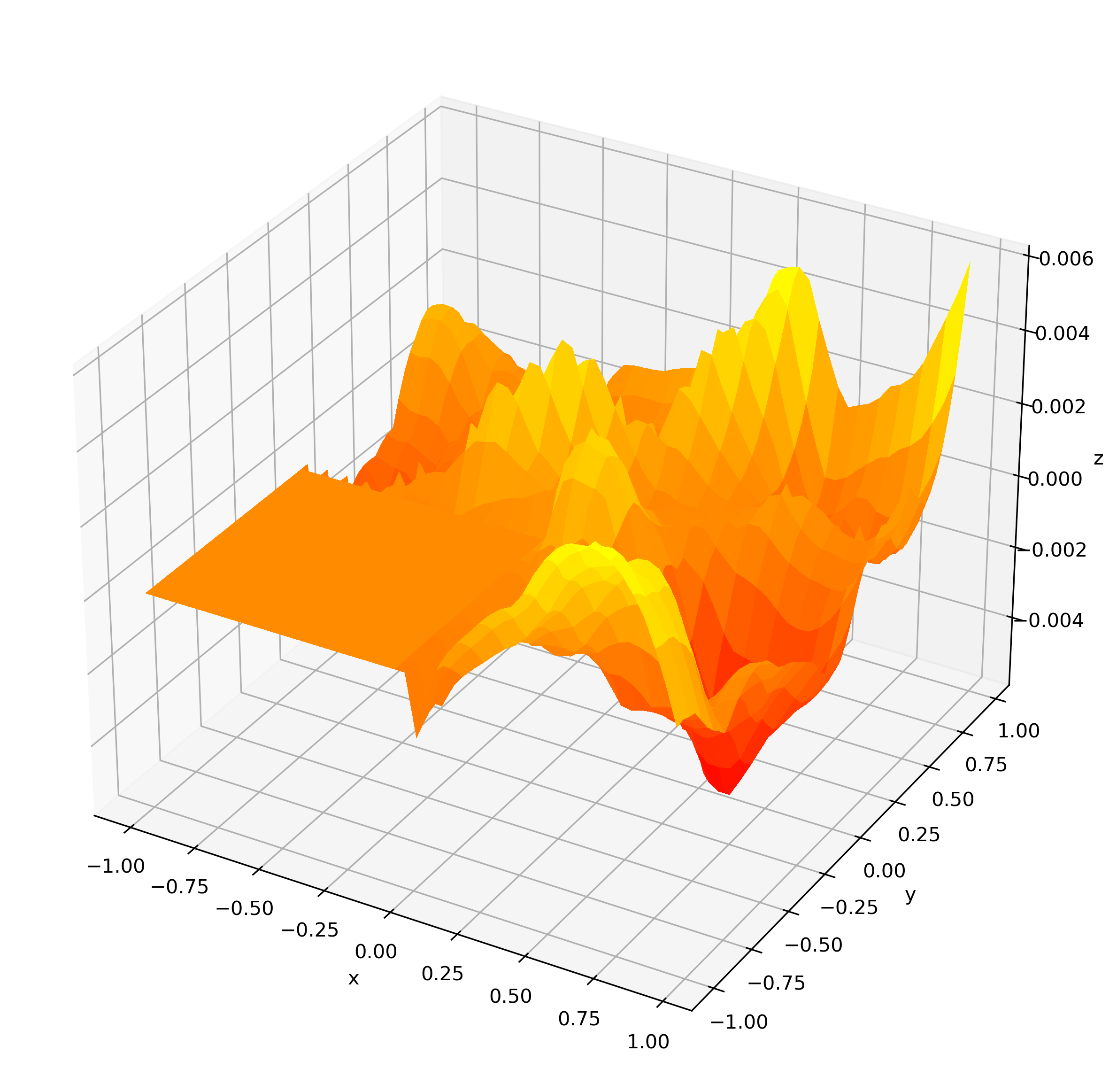}
    \caption{Visualization of the difference between the solution obtained with classical solution and the solution from the PINN. The mean squared error is approximately $2.43\text{e-}6$.}
    \label{fig:heat-error-pinn}
\end{figure}

A physics-informed neural network has been trained to directly approximate the result, respecting the heat transport differential equation, the heat source function, and boundary conditions presented in subsection \ref{sub:coefficient-approximation}. %Figure \ref{fig:heat-pinn} shows the solution obtained from the classical method in the blue-green palette. The solution obtained from the PINN is overlaid in red-yellow colors. 
A difference between the solution obtained with the IGA solver \ref{IGA} and the solution obtained from the PINN is shown in Figure \ref{fig:heat-error-pinn}. The results are satisfactory, unfortunately PINNs can't learn families of functions, at least not in their original form. 

\subsection{A short note on overfitting and network capacity}
Despite the loss function plots suggesting the overfitting of models, all networks performed well on never–seen-before data with MSE values of the same order of magnitude as in the training data.
Coefficient approximating models were also evaluated with a 5-fold cross validation. 
While the network architectures might suggest reduced network capacity compared to the architectures with the same number of input nodes as the number of input parameters, adding dropout layers hasn't improved the accuracy. 

%\begin{figure}[H]
%    \centering
%    \includegraphics[width=0.8\linewidth]{pinn_approximation}
%    \caption{Visualization of the heat equation solution. The solution obtained from the classical method in the blue-green palette, while the solution obtained from the PINN is overlaid in red-yellow colors. The mean squared error is approximately $2.43\text{e-}6$.}
%    \label{fig:heat-pinn}
%\end{figure}

\begin{figure}[H]
    \centering
    \includegraphics[width=\linewidth]{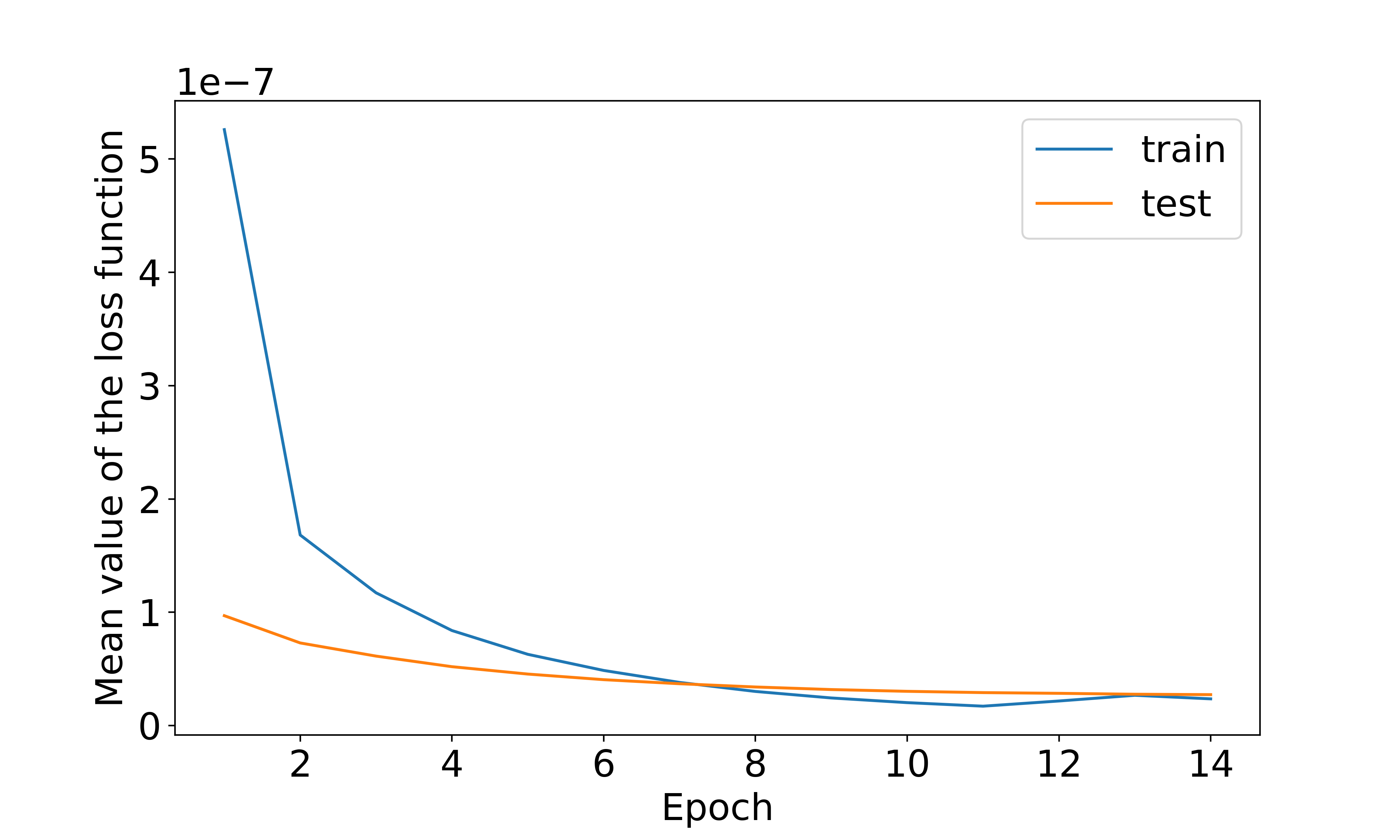}
    \caption{Mean squared error loss function over the 15 training epochs. Each batch contained 40 entries. }
    \label{fig:heat-pinn-direct-loss}
\end{figure}

\section{Conclusions and future work}

In terms of research objectives, the results of this paper are as follows:

\begin{enumerate}
    \item \textbf{Investigating the viability of using NN for approximating coefficients of base B-spline functions used in IGA.}
    
    Application of neural networks in approximating coefficients of the base functions proved to be possible. Deep Neural Networks were able to learn a family of functions with more than acceptable accuracy for CAD applications such as prototyping. More research is needed to check whether neural networks can accurately approximate coefficients with arbitrary knot vectors as an input, as opposed to a predefined knot vectors used in this paper.
    
    \item \textbf{Comparing the two approaches in applying NNs to IGA in terms of accuracy and training time.}
    
    Approximating base functions' coefficients was a more accurate (figures \ref{fig:fem-dnn-solution}, \ref{fig:error-smooth-PINN}, \ref{fig:error-PINN}, %\ref{fig:heat-pinn} 
and \ref{fig:heat-error-pinn}) and quicker (figure \ref{fig:heat-dnn-coeff-loss}, \ref{fig:heat-dnn-direct-loss} and \ref{fig:heat-pinn-direct-loss}) to train approach, but the network has to be retrained for different meshes. Direct approximation with deep neural networks was slower to train (table \ref{tab:heat-times-comparison}) but was independent of the underlying mesh, as the networks approximate the final solution. Physics-informed Neural Networks could only predict one specific function, but it required less data from solved problems, provided all the underlying equations and boundary conditions are known at training time. 
    
\end{enumerate}

Further research is needed in some areas. The question whether the number of B-splines and coefficients along with the size of the network impact the approximation error is left unanswered. The data summarized in the table \ref{tab:heat-dnn-coeff-results} suggest that there might exist a correlation between the problem size, the desired accuracy and the approximation error, due to the more difficult learning process which involves more variables. 

\section*{Acknowledgments}
The paper was partially financed by AGH University of Science and Technology Statutory Fund.

\section{Appendix A: One-dimensional example of neural network learning coefficients of B-splines}

\subsection{One dimensional heat-transfer problem}
Let us introduce the knot vector [0 0 0 1 1 1] defining the quadratic B-spline basis functions with $C^0$ separators
\begin{equation}
B_{1,2}(x)=(1-x)^2; \quad B_{2,2}(x)=2x(1-x); \quad B_{3,2}(x)=x^2
\end{equation}

\begin{figure}[h]
\includegraphics[scale=0.4]{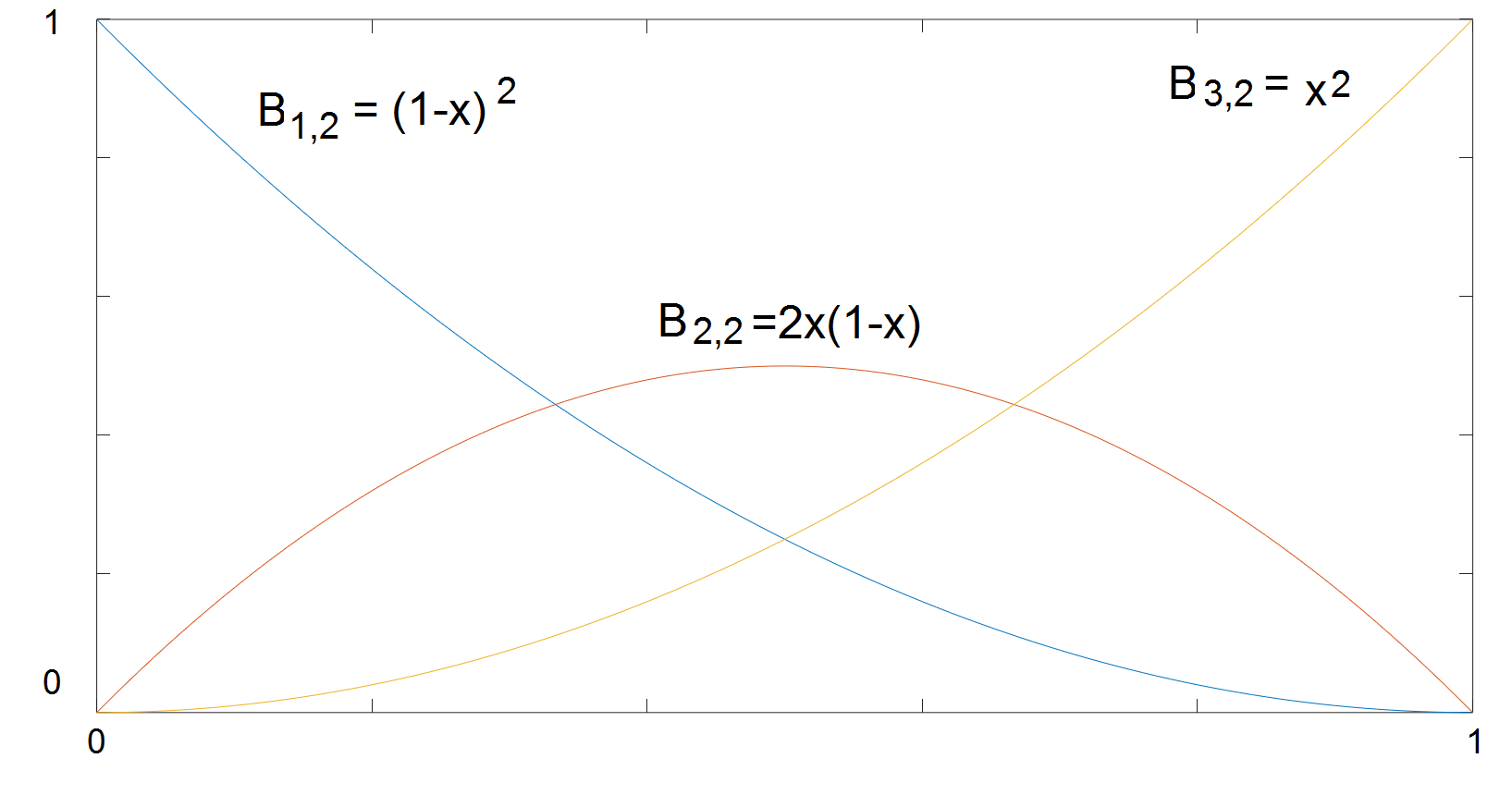}
\caption{Three B-splines over a single interval (element)}
\end{figure}

Let us introduce the problem
\begin{equation}
-u''(x)=f(x) \quad x \in (0,0.5)
\end{equation}
defined over $x \in (0,0.5)$,
with boundary conditions $u(0)=0$ and $u'(0.5)=g(x)$. 
We setup $g(x)=n \pi cos(n \pi x)$ and $f(x)=n^2 \pi^2 sin(n \pi x)$.
The family of solution of this problem are
\begin{equation}f_n(x)=sin(n \pi x) \end{equation}
\begin{figure}
\includegraphics[scale=0.4]{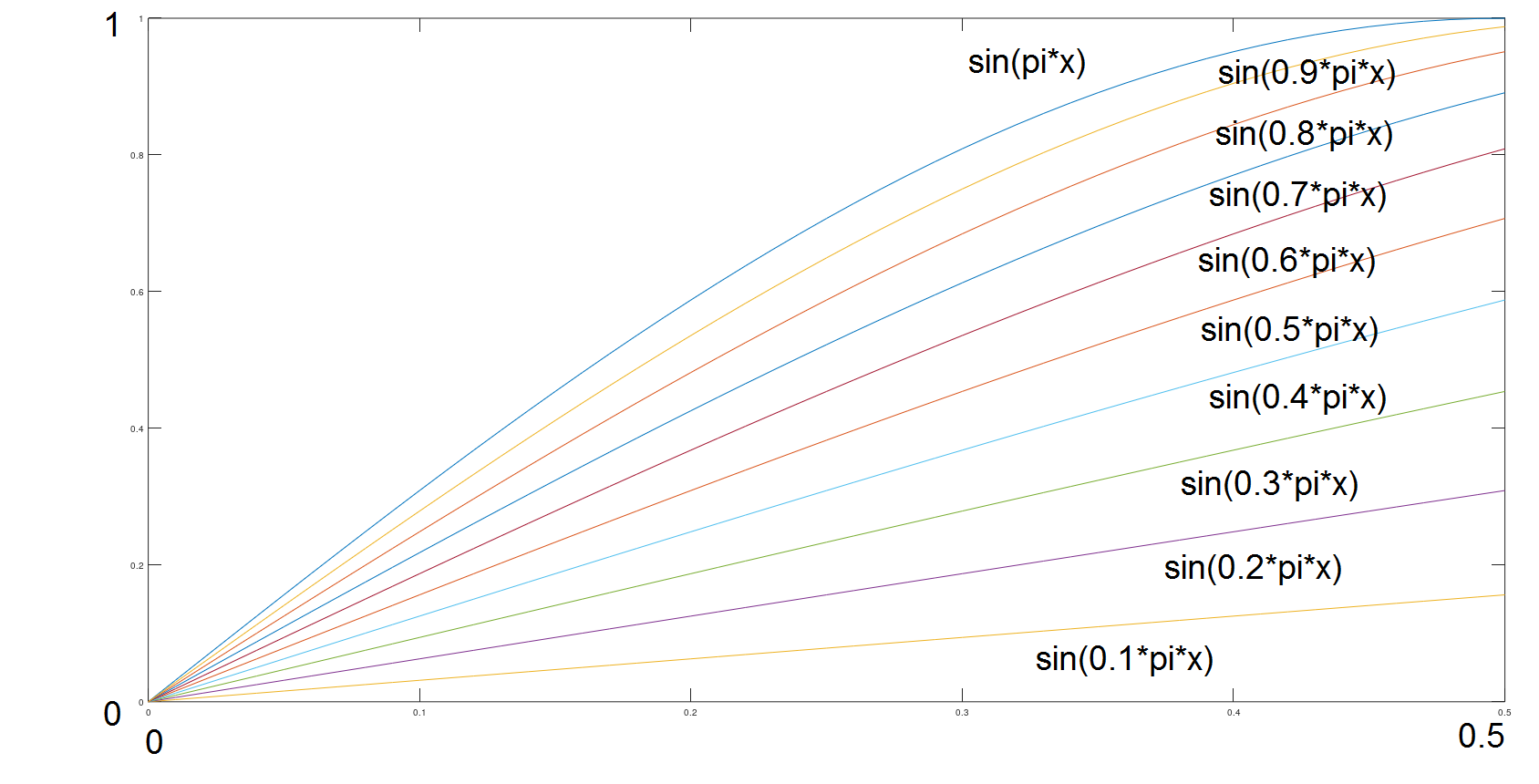}
\caption{Plot of different solutions $f_n(x)$ for $n=0.1,0.2,...,1$.}
\end{figure}

We transform this problem into the weak form
\begin{eqnarray}
\int_0^{0.5} u'(x)v'(x) dx = \int_0^{0.5} f(x) v(x) dx+v(0.5)g(0.5) \quad \forall v
\end{eqnarray}
and we discretize with B-spline basis functions
\begin{eqnarray}
u_h=\sum_{i=1,2,3} u_i B_{i,2}(x) 
\end{eqnarray}
to obtain
\begin{eqnarray}
\begin{bmatrix} \int_{0,1/2} B_{1,2}'(x)B_{1,2}'(x) dx & \int_{0,1/2} B_{1,2}'(x)B_{2,2}'(x) dx  & \int_{0,1/2} B_{1,2}'(x)B_{3,2}'(x) dx \\
\int_{0,1/2} B_{2,2}'(x)B_{1,2}'(x) dx & \int_{0,1/2} B_{2,2}'(x)B_{2,2}'(x) dx  & \int_{0,1/2} B_{2,2}'(x)B_{3,2}'(x) dx \\
\int_{0,1/2} B_{3,2}'(x)B_{1,2}'(x) dx & \int_{0,1/2} B_{3,2}'(x)B_{2,2}'(x) dx  & \int_{0,1/2} B_{3,2}'(x)B_{3,2}'(x) dx \\
\end{bmatrix} \begin{bmatrix} u_1 \\ u_2 \\ u_3 \end{bmatrix} = \nonumber \\ \begin{bmatrix} \int_{0,1/2} B_{1,2}(x) f_n(x) dx  \\ 
\int_{0,1/2} B_{2,2}(x) f_n(x) dx \\
\int_{0,1/2} B_{3,2}(x) f_n(x) dx+n\pi cos(n\pi 0.5) \end{bmatrix} 
\label{eq:problem}
\end{eqnarray}
%\begin{eqnarray}
%\begin{bmatrix} \frac{1}{5} & \frac{1}{10} & \frac{1}{30} \\
%\frac{1}{10} & \frac{2}{15} & \frac{1}{10} \\
%\frac{1}{30} & \frac{1}{10} & \frac{1}{5} 
%\end{bmatrix}  \begin{bmatrix} u_1 \\ u_2 \\ u_3 \end{bmatrix} = \begin{bmatrix} \int_{0,1/2} (1-x^2) sin(n\Pi *x) dx  \\ 
%\int_{0,1/2} 2x(1-x) sin(n\Pi *x) dx \\
%\int_{0,1/2} x^2 sin(n\Pi *x) dx \end{bmatrix} 
%\end{eqnarray}
%\begin{eqnarray}
%\begin{bmatrix} \frac{1}{5} & \frac{1}{10} & \frac{1}{30} \\
%\frac{1}{10} & \frac{2}{15} & \frac{1}{10} \\
%\frac{1}{30} & \frac{1}{10} & \frac{1}{5} 
%\end{bmatrix} \begin{bmatrix} u_1 \\ u_2 \\ u_3 \end{bmatrix} = \begin{bmatrix} \frac{\Pi^2 n^2+2cos (\Pi n)-2}{\Pi^3 n^3} \\ 
%\frac{-2\Pi n sin (\Pi n)-4 cos (\Pi n)+4}{\Pi^3 n^3} \\
%\frac{(2-\Pi^2 n^2) cos (\Pi n)+2\Pi n sin(\Pi n)-2}{\Pi^3 n^3} \\
%\end{bmatrix} 
%\end{eqnarray}
\begin{figure}[h]
\includegraphics[scale=0.4]{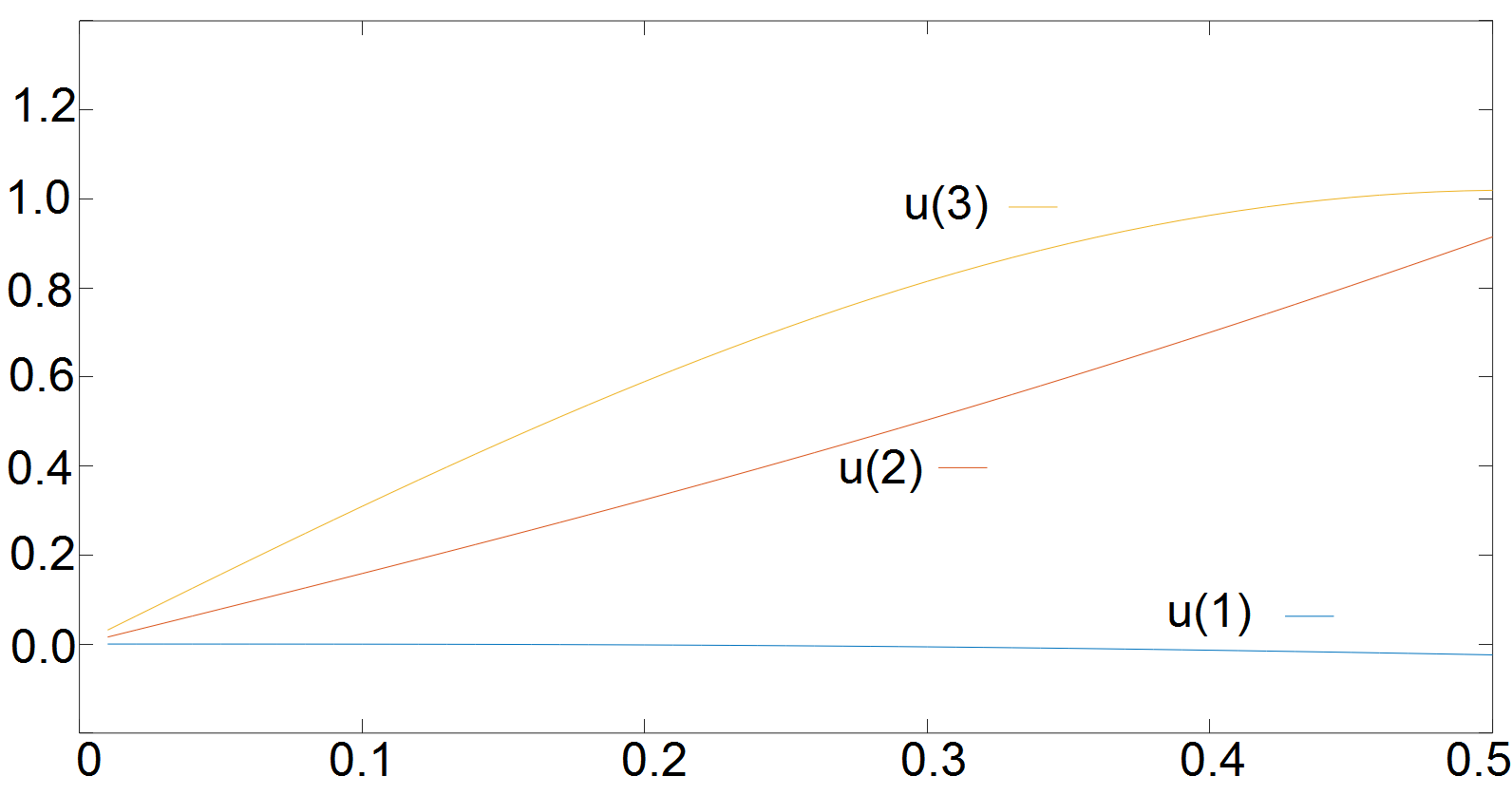}
\caption{Coefficients of approximation $u_1B_{1,2}(x)+u_2B_{2,2}+u_3B_{3,2}$ for $n\in (0,0.5)\subset {\cal R}$.}
\end{figure}

\subsection{Artificial neural network for $u_i$}

Let us introduce the artificial neural network
\begin{equation}
ANN_i(n)=u_i \label{ANN1}
\end{equation}
where $n$ is the index of the $f_n$ function, $i=1,2,3$ (for three coefficients of B-splines).

Given sinus family function index $n$, it returns the coefficient $u_i$ of B-splines for approximation of this function over $(0,0.5)$.

\begin{equation}
ANN_i(n)=c_i \sigma \left( a_i n +b_i \right)+d_i
\end{equation}
where the activation function
\begin{equation}
\sigma(x)=\frac{1}{1+e^{-x}}
\end{equation}

\subsection{Training}

The goal of the training is to find values of the weights $a_i,b_i,c_i,d_i$

We prepare a set of samples

\begin{itemize}
\item We randomly select $n \in (0,1)$
\item We solve the IGA problem (\ref{eq:problem})
\item Input data $(n)$, output data $(u_1,u_2,u_3)$
\end{itemize}
In other words, we train the neural network for some selected functions $f_n$, with a hope, that it will work for a given function of interest from the family.

How to train the artificial neural network?
We define the error function
\begin{eqnarray}
e_i(n)=0.5\left(ANN_i(n)-u_i(n)\right)^2=0.5\left(c_i \sigma \left( a_i n +b_i \right)+d_i-u_i(n)\right)^2= \\
0.5\left(\left(\frac{c_i} {1+exp(-\left( a_i n +b_i \right)}+d_i\right)-u_i(n)\right)^2
\end{eqnarray}
Now, we compute the derivatives
\begin{eqnarray}
\frac{\partial e_i(n)}{\partial a_i} = \frac{c_i n exp(-a_in-b_i)(ANN_i(n)-u_i(n))}{(exp(-ax-b)+1)^2} \\
\frac{\partial e_i(n)}{\partial b_i} = \frac{c_i exp(-a_in-b_i)(ANN_i(n)-u_i(n))}{(exp(-ax-b)+1)^2} \\
\frac{\partial e_i(n)}{\partial c_i} = \frac{(ANN_i(n)-u_i(n))}{(exp(-ax-b)+1)} \\
\frac{\partial e_i(n)}{\partial d_i} = (ANN_i(n)-u_i(n)) \\
\end{eqnarray}
they say ``how fast the error is changing if I modify a given coefficient''.

We loop through the data set $\{ n, (u_1(n),u_2(n),u_3(n))\}_{n \in A}$ where $A$ is the set of selected points from $(0,0.5)$, and we train each of the three $ANN_1$, $ANN_2$, and $ANN_3$

\begin{enumerate}
\item Select $(n,(u_1,u_2,u_3))$
\item Compute $u_i=ANN_i(n)=c_i \sigma \left( a_i n +b_i \right)+d_i$
\item Compute  $e_i(n)$
\item Compute  $\frac{\partial e_i(n)}{\partial a_i}, \frac{\partial e_i(n)}{\partial b_i}, \frac{\partial e_i(n)}{\partial c_i}, \frac{\partial e_i(n)}{\partial d_i}$
\item Correct 
\begin{eqnarray}
a_i = a_i - \eta * \frac{\partial e_i(n)}{\partial a_i} \\
b_i = b_i - \eta * \frac{\partial e_i(n)}{\partial b_i} \\
c_i = c_i - \eta * \frac{\partial e_i(n)}{\partial c_i} \\
d_i = d_i - \eta * \frac{\partial e_i(n)}{\partial d_i} 
\end{eqnarray}
\end{enumerate}
where $\eta \in (0,1)$. This is like a local gradient method.

\subsection{MATLAB implementation}

      {\tt \inred{\% Creation of dataset}} \\
      {\tt A = [1/5 1/10 1/30; 1/10 2/15 1/10;  1/30 1/10 1/5];}\\
      {\tt i=1;}\\
      {\tt for n=0.01:0.01:0.5}\\
      {\tt rhs= [  (pi*pi*n*n+2*cos(pi*n)-2)/(pi*pi*pi*n*n*n);}\\
      {\tt (-2*pi*n*sin(pi*n)-4*cos(pi*n)+4)/(pi*pi*pi*n*n*n); }\\
      {\tt ((2-pi*pi*n*n)*cos(pi*n)+2*pi*n*sin(pi*n)-2)/(pi*pi*pi*n*n*n) ];}\\
      {\tt u=A $\backslash$ rhs;}\\
      {\tt dataset\_in(i)=n; }\\
      {\tt dataset\_u1(i)=u(1); }\\
      {\tt dataset\_u2(i)=u(2); }\\
      {\tt dataset\_u3(i)=u(3); }\\
      {\tt i=i+1;}\\
      {\tt endfor}\\
      {\tt ndataset=i-1;}

      {\tt \inred{\% Training}} \\
      {\tt a1=1.0; b1=1.0; c1=1.0; d1=1.0;}\\
      {\tt eta1=0.1;}\\
      {\tt r = 0 + (1-0).*rand(ndataset,1);}\\
      {\tt r=r.*ndataset;}\\
      {\tt for j=1:ndataset}\\
      {\tt   i=floor(r(j));}\\
      {\tt   eval1 = c1*1.0/(1.0+exp(-(a1*dataset\_in(i)+b1)))+d1;}\\
      {\tt   error1 = 0.5*(eval1-dataset\_u1(i))$^2$;}\\
      {\tt   derrorda = c1*dataset\_in(i)*exp(-a1*dataset\_in(i)-b1)*}\\
      {\tt \; \; \; (eval1-dataset\_u1(i))/(exp(-a1*dataset\_in(i)-b1)+1)$^2$;}\\
      {\tt   a1=a1-eta1*  derrorda;}\\
      {\tt   derrordb = c1*exp(-a1*dataset\_in(i)-b1)*}\\
      {\tt \; \; \; (eval1-dataset\_u1(i))/(exp(-a1*dataset\_in(i)-b1);}
      {\tt   b1=b1-eta1*  derrordb;}\\
      {\tt   derrordc = (eval1-dataset\_u1(i))/(exp(-a1*dataset\_in(i)-b1)+1);}\\
      {\tt   c1=c1-eta1*  derrordc;}\\
      {\tt   derrordd = (eval1-dataset\_u1(i));}\\
      {\tt   d1=d1-eta1*  derrordd;}\\

      {\tt   \inred{\% evaluation of ANN approximation of sin(n*pi*x) for n=0.333}}\\
      {\tt   n=0.333;}\\
      {\tt   u1 = c1*1.0/(1.0+exp(-(a1*n+b1)))+d1;}\\
      {\tt   u2 = c2*1.0/(1.0+exp(-(a2*n+b2)))+d2;}\\
      {\tt   u3 = c3*1.0/(1.0+exp(-(a3*n+b3)))+d3;}\\
      {\tt   x=0:0.01:0.5;}\\
      {\tt   y=sin(n*pi.*x);}\\
      {\tt   z=u1*(1-x).$^2$+u2*2*x.*(1-x)+u3*x.$^2$;}\\
      {\tt   plot(x,y,x,z);}\\

We tried starting points $1.0,10.0,-1.0,-10.0$ for all the combinations of $a_i,b_i,c_i,d_i$ (256 runs) and the best result (smaller errors) we obtain for 

$a_1=b_1=c_1=d_1=1.0$; 
$a_2=b_2=1,c_2=10.0,d_2=-1.0$
$a_3=b_2=3,c_3=10.0,d_3=-1.0$

We used $\eta=0.1$. We coded the ANN and the training in hand-made MATLAB code.

\begin{figure}[h]%
\includegraphics[scale=0.4]{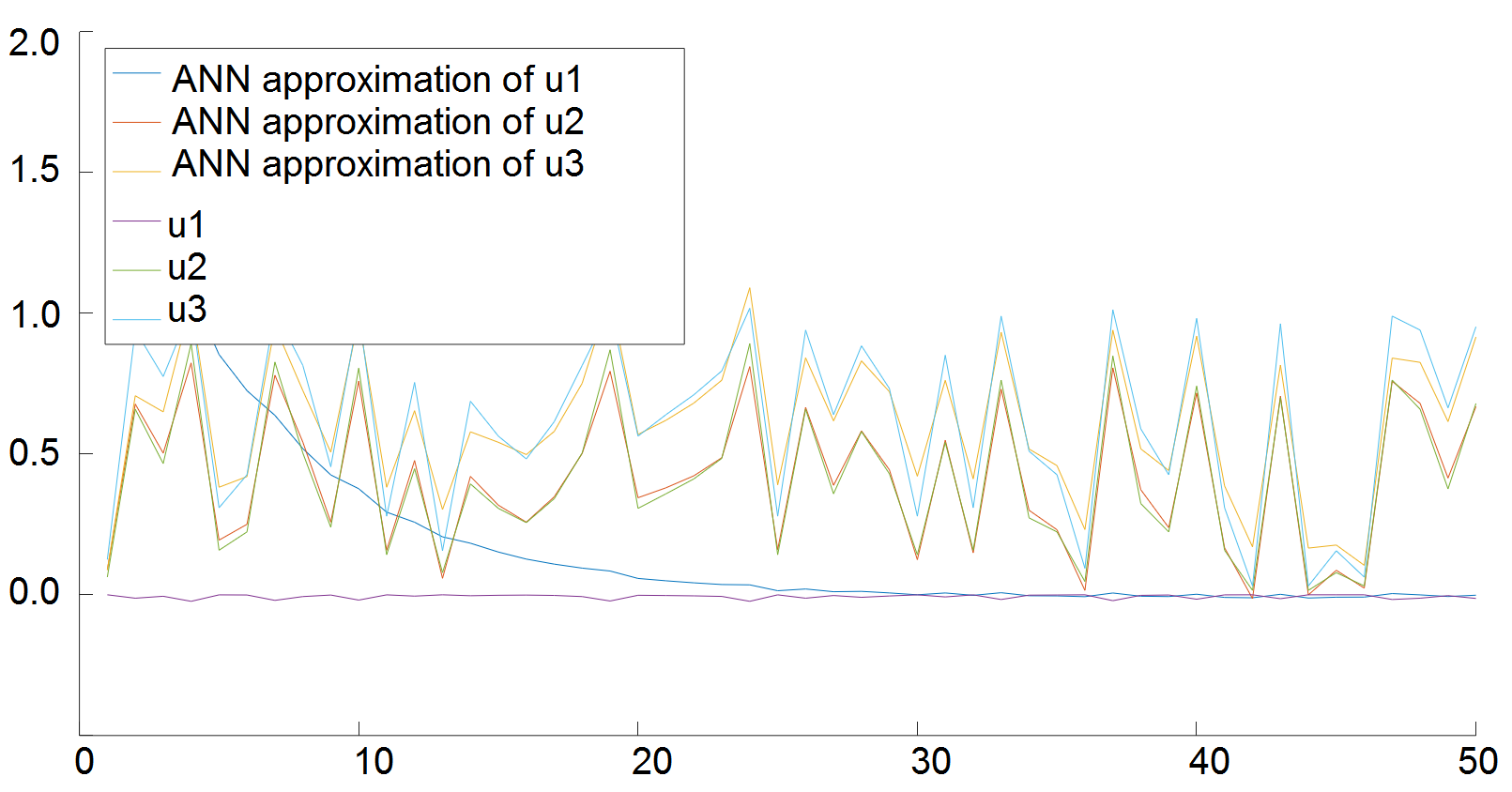}
\caption{Training for the simple artificial nerual network (\ref{ANN1})
 starting from $a_1=b_1=c_1=d_1=1.0$; $a_2=b_2=1,c_2=10.0,d_2=-1.0$, $a_3=b_2=3,c_3=10.0,d_3=-1.0$, for $\eta=0.1$.} \label{t1}
\end{figure}
\begin{figure}[h]
\includegraphics[scale=0.4]{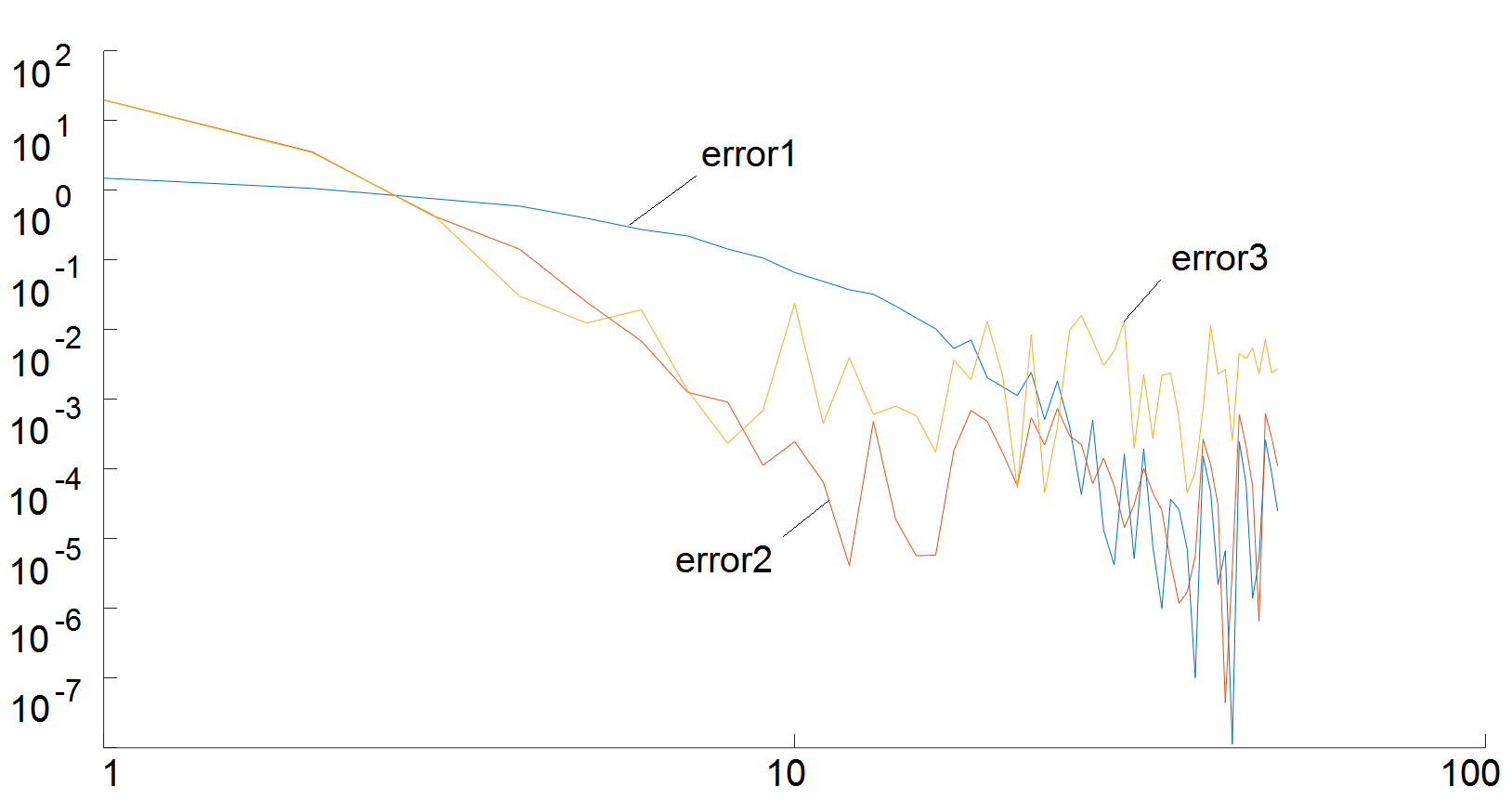}
\caption{Convergence of errors for training of ANN1, ANN2, ANN3} \label{t2}
% starting from 1,1,10,-10}
\end{figure}
%\begin{figure}[h]
%\includegraphics[scale=0.4]{min3.png}
%\caption{Sensitivity analysis for $a$ and $b$ parameters.} \label{sens1}
%\end{figure}
%\begin{figure}[h]
%\includegraphics[scale=0.4]{min2.png}
%\caption{Sensitivity analysis  for $c$ and $d$ parameters.} \label{sens2}
%\end{figure}

\subsection{Verification}

\begin{figure}[h]
\includegraphics[scale=0.4]{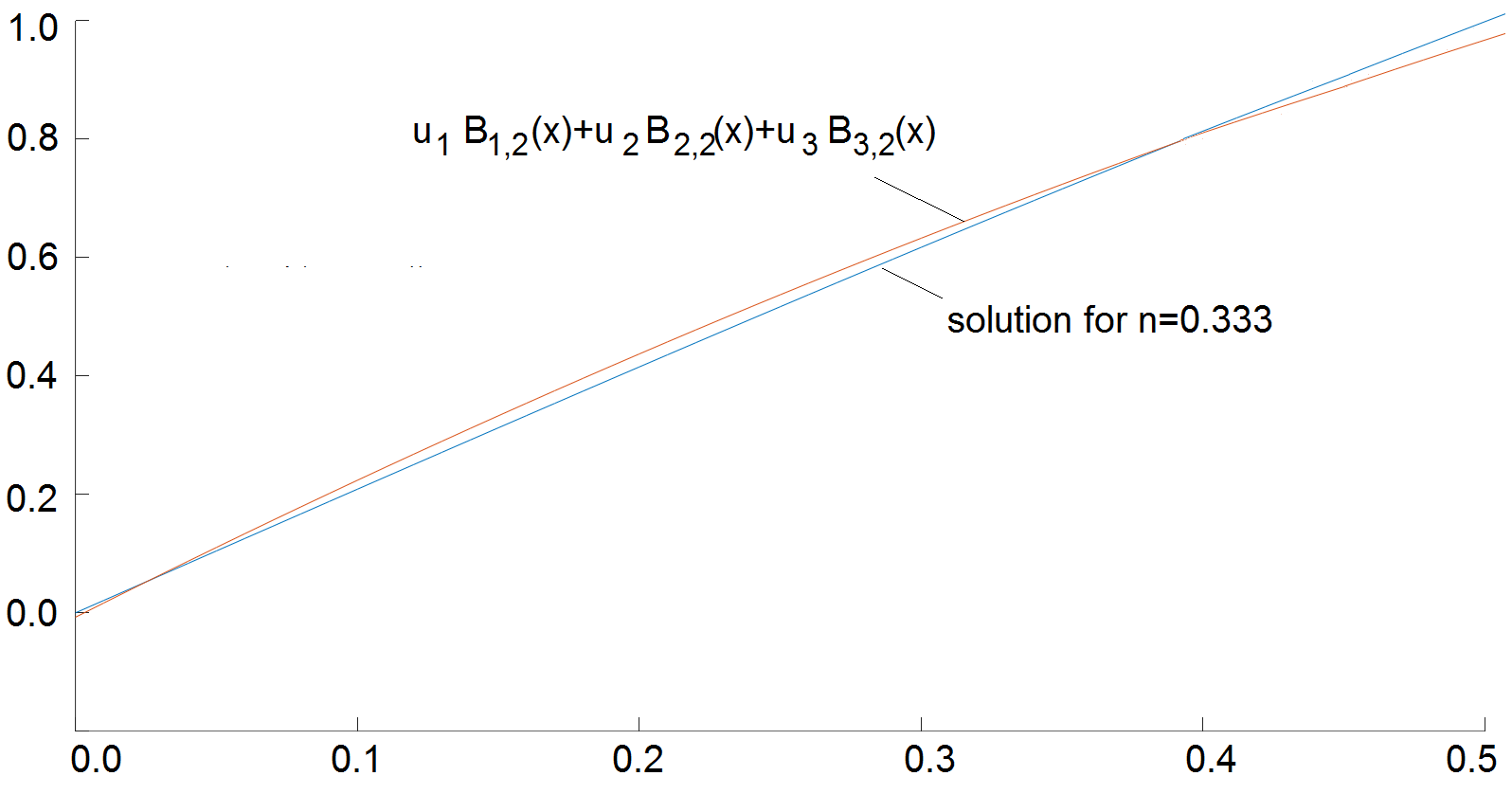}
\caption{Verification of the neural network approximation of solution for $n=0.333$ with
      {\tt  u1 = c1*1.0/(1.0+exp(-(a1*n+b1)))+d1; u2 = c2*1.0/(1.0+exp(-(a2*n+b2)))+d2;}
      {\tt   u3 = c3*1.0/(1.0+exp(-(a3*n+b3)))+d3; z=u1*(1-x).$^2$+u2*2*x.*(1-x)+u3*x.$^2$;}}\label{0333}
\end{figure}

In Figure \ref{t1} we present the training over 50 samples, and in Figure \ref{t2} we present the convergence of the training.

%In Figure \ref{sens1}-\ref{sens2} we plot the error functions for parameters $a,b$ and $c,d$ for the second B-spline coefficients.

We select $n=0.333$ and we compute

\begin{eqnarray}
z(x)=ANN_1(n)*B_{1,2}(x)+ANN_2(n)*B_{2,2}(x)+ANN_3(n)*B_{3,2}(x) = \\
\left(\frac{c_1}{1+exp(-a_1*n-b_1)}+d_1\right)*(1-x)^2+ \left(\frac{c_2}{1+exp(-a_2*x-b_2)}+d_2\right)*2x(1-x)^2+\\ +\left(\frac{c_3}{1+exp(-a_3*x-b_3)}+d_3\right)*x^2 
\end{eqnarray}

we compare with $sin(0.333\pi x)$ and $sin(0.777\pi x)$ in Figure \ref{0333}.

%\begin{figure}[h]
%\includegraphics[scale=0.4]{error_ANN.png}
%\caption{Error between the Artificial neural network and L2 projection $\| ANN(i) B_{i,2}(x) -  \sum_{i=1,2,3} u_i B_{i,2}(x)\|^2_{L^2} = \| ANN-IGA \|$, as well as between the Artificial Neural Network and exact solution $\| sin(n*\pi*x)- ANN(i) B_{i,2}(x)\|^2_{L^2}=\|exact-ANN\|$, for different values of parameter $n$ in $sin(n*\pi*x)$.}
%\end{figure}

\section{Appendix B: One-dimensional example of neural network learning IGA solution}

\subsection{One dimensional heat-transfer problem}

We focus again on the heat-transfer problem
\begin{equation}
-u''(x)=f(x) \quad x \in (0,0.5)
\end{equation}
defined over $x \in (0,0.5)$,
with boundary conditions $u(0)=0$ and $u'(0.5)=g(x)$, with $g(x)=n \pi cos(n \pi x)$ and $f(x)=n^2 \pi^2 sin(n \pi x)$.
We transform this problem into the weak form and we discretize with B-spline basis functions
\begin{eqnarray}
\begin{bmatrix} \int_{0,1/2} B_{1,2}'(x)B_{1,2}'(x) dx & \int_{0,1/2} B_{1,2}'(x)B_{2,2}'(x) dx  & \int_{0,1/2} B_{1,2}'(x)B_{3,2}'(x) dx \\
\int_{0,1/2} B_{2,2}'(x)B_{1,2}'(x) dx & \int_{0,1/2} B_{2,2}'(x)B_{2,2}'(x) dx  & \int_{0,1/2} B_{2,2}'(x)B_{3,2}'(x) dx \\
\int_{0,1/2} B_{3,2}'(x)B_{1,2}'(x) dx & \int_{0,1/2} B_{3,2}'(x)B_{2,2}'(x) dx  & \int_{0,1/2} B_{3,2}'(x)B_{3,2}'(x) dx \\
\end{bmatrix} \begin{bmatrix} u_1 \\ u_2 \\ u_3 \end{bmatrix} = \nonumber \\ \begin{bmatrix} \int_{0,1/2} B_{1,2}(x) f_n(x) dx  \\ 
\int_{0,1/2} B_{2,2}(x) f_n(x) dx \\
\int_{0,1/2} B_{3,2}(x) f_n(x) dx+n\pi cos(n\pi 0.5) \end{bmatrix} 
\label{eq:problem}
\end{eqnarray}

\subsection{Artificial neural network approximating solution}

Let us introduce the artificial neural network
\begin{equation}
ANN(n,x)=y \label{ANN2}
\end{equation}
where $n$ is the index of the $f_n$ function, and $x$ is the argument.

\begin{equation}
ANN(n,x)=c \sigma \left( \begin{bmatrix} a_1 a_2 \end{bmatrix} \begin{bmatrix} n \\ x \end{bmatrix} +b \right)+d
\end{equation}
where the activation function
\begin{equation}
\sigma(x)=\frac{1}{1+e^{-x}}
\end{equation}

\subsection{Training}

The goal of the training is to find values of the weights $a_1,a_2,b,c,d$

We prepare a set of samples

\begin{itemize}
\item We randomly select $n \in (0,1)$ and $x \in (0,1)$
\item We solve the IGA problem (\ref{eq:problem}) to obtain $\left(u_1,u_2,u_3\right)$
\item Input data $(n,x)$, output data $y(n,x)=u_1 B_{1,2}+u_2B_{2,2}+u_3B_{3,2}$
\end{itemize}

We define the error function
\begin{eqnarray}
e(n,x)=0.5\left(ANN(n,x)-y(n,x)\right)^2=0.5\left(c \sigma \left( a_1 n +a_2x+ b \right)+d-y(n,x)\right)^2= \nonumber \\
0.5\left(\left(\frac{c} {1+exp(-a_1 n - a_2 x - b )}+d\right)-y(n,x)\right)^2
\end{eqnarray}
Now, we compute the derivatives
\begin{eqnarray}
\frac{\partial e(n,x)}{\partial a_1} = \frac{cn exp \left(-a_1n-a_2x-b\right)\left(\frac{c}{exp\left(-a_1n-a_2x-b\right)+1}+d-y\right)}{\left(exp\left(-a_1n-a_2x-b\right)+1\right)^2}
 \\
\frac{\partial e(n,x)}{\partial a_2} = \frac{cx exp \left(-a_1n-a_2x-b\right)\left(\frac{c}{exp\left(-a_1n-a_2x-b\right)+1}+d-y\right)}{\left(exp\left(-a_1n-a_2x-b\right)+1\right)^2} \\
\frac{\partial e(n,x)}{\partial b} = \frac{c exp \left(-a_1n-a_2x-b\right)\left(\frac{c}{exp\left(-a_1n-a_2x-b\right)+1}+d-y\right)}{\left(exp\left(-a_1n-a_2x-b\right)+1\right)^2} \\
\frac{\partial e(n,x)}{\partial c} = \frac{\frac{c}{exp\left(-a_1n-a_2x-b\right)+1}+d-y}{exp\left(-a_1n-a_2x-b\right)+1} \\
\frac{\partial e(n,x)}{\partial d} = \frac{c}{exp\left(-a_1n-a_2x-b\right)+1}+d-y
\end{eqnarray}

We loop through the data set $\{ (n,x), y\}_{(n,x) \in A}$ where $A$ is the set of selected points from $(0,0.5)\times(0,1)$, and we train $ANN$

\begin{enumerate}
\item Select $((n,x),y)$
\item Compute $y=ANN(n,x)=c \sigma \left( \begin{bmatrix} a_1 a_2 \end{bmatrix} \begin{bmatrix} n \\ x \end{bmatrix} +b \right)+d$
\item Compute  $e(n,x)$
\item Compute  $\frac{\partial e(n,x)}{\partial a_1}, \frac{\partial e(n,x)}{\partial a_2}, \frac{\partial e(n,x)}{\partial b}, \frac{\partial e(n,x)}{\partial c}, \frac{\partial e(n,x)}{\partial d}$
\item Correct 
\begin{eqnarray}
a_1 = a_1 - \eta * \frac{\partial e(n,x)}{\partial a_1} \\
a_2 = a_2 - \eta * \frac{\partial e(n,x)}{\partial a_2} \\
b = b - \eta * \frac{\partial e(n,x)}{\partial b} \\
c = c - \eta * \frac{\partial e(n,x)}{\partial c} \\
d = d - \eta * \frac{\partial e(n,x)}{\partial d} 
\end{eqnarray}
\end{enumerate}
where $\eta \in (0,1)$. 

\begin{figure}[h]%
\includegraphics[scale=0.4]{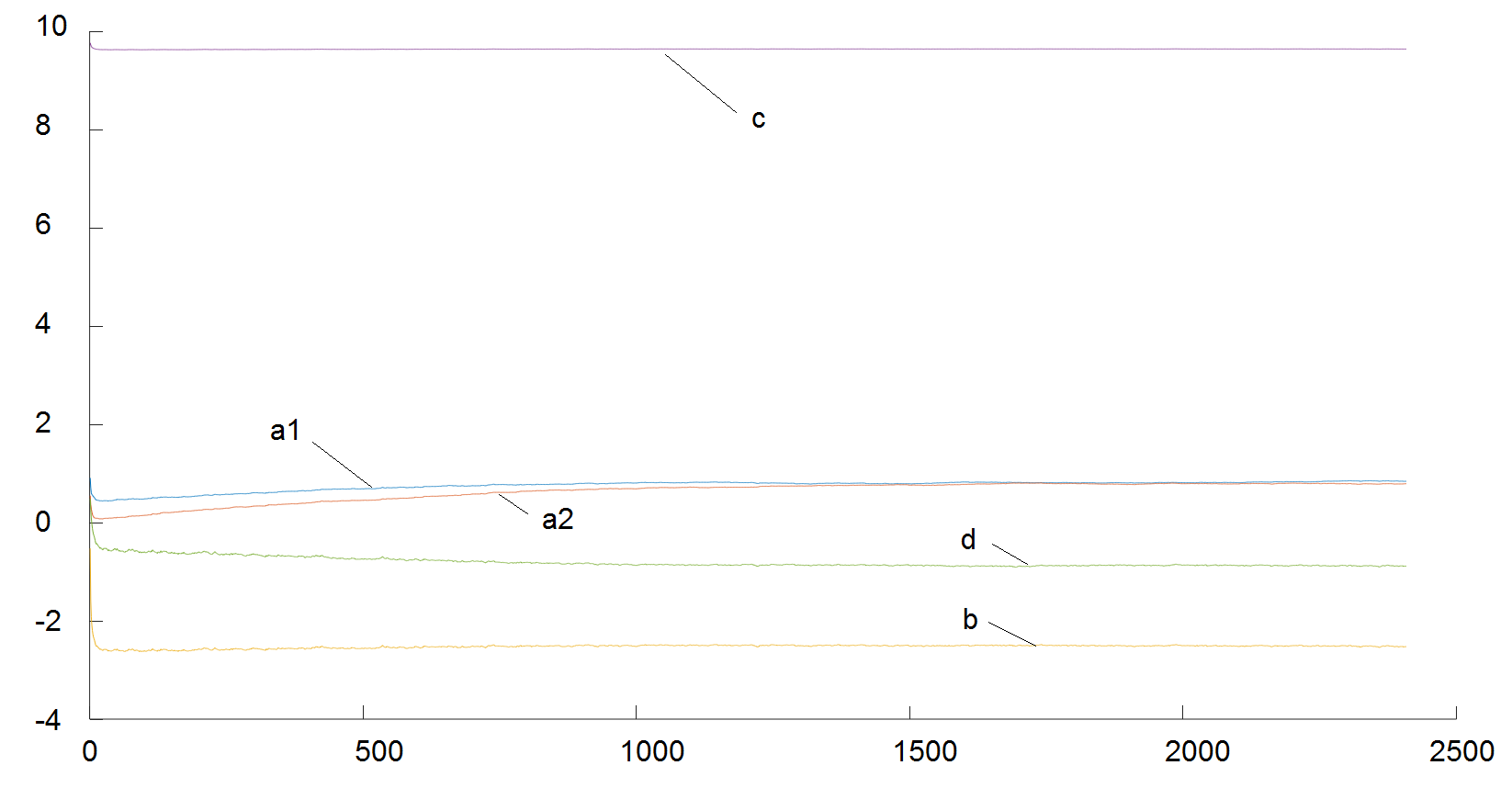}
\caption{Training for the simple artificial neuRal network (\ref{ANN2})
 starting from $a_1=a_2=b=c=d=1.0$, for $\eta=0.1$.} \label{t1a}
\end{figure}
\begin{figure}[h]
\includegraphics[scale=0.4]{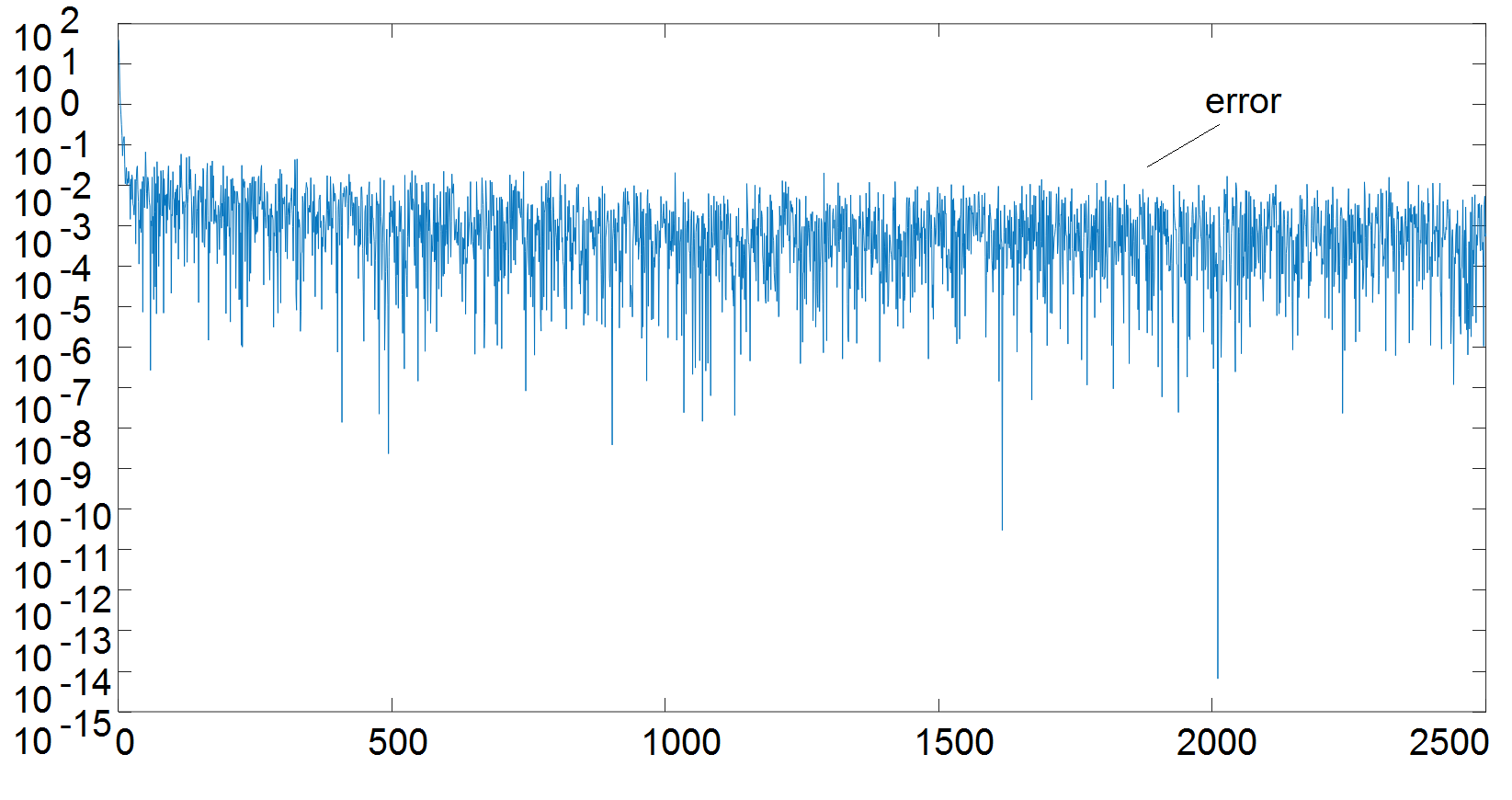}
\caption{Convergence of error for training of ANN} \label{t2a}
% starting from 1,1,10,-10}
\end{figure}

\subsection{MATLAB implementation}

      {\tt \inred{\% Creation of dataset}} \\
      {\tt A = [1/5 1/10 1/30; 1/10 2/15 1/10;  1/30 1/10 1/5];}\\
      {\tt i=1;}\\
      {\tt for n=0.01:0.001:0.5}\\
      {\tt for x=0.01:0.001:0.5}\\
      {\tt rhs= [  (pi*pi*n*n+2*cos(pi*n)-2)/(pi*pi*pi*n*n*n);}\\
      {\tt (-2*pi*n*sin(pi*n)-4*cos(pi*n)+4)/(pi*pi*pi*n*n*n); }\\
      {\tt ((2-pi*pi*n*n)*cos(pi*n)+2*pi*n*sin(pi*n)-2)/(pi*pi*pi*n*n*n) ];}\\
      {\tt u=A $\backslash$ rhs;}\\
      {\tt y=u(1)*(1-x).$^2$+u(2)*2*x.*(1-x)+u(3)*x.$^2$;}\\
      {\tt dataset\_in\_n(i)=n; }\\
      {\tt dataset\_in\_x(i)=x; }\\
      {\tt dataset\_y(i)=y; }\\
      {\tt i=i+1;}\\
      {\tt endfor}\\
      {\tt endfor}\\
      {\tt ndataset=i-1;}

      {\tt \inred{\% Training}} \\
      {\tt a1=1.0; a2=1.0; b=1.0; c=1.0; d=1.0;}\\
      {\tt eta=0.1;}\\
      {\tt r = 0 + (1-0).*rand(ndataset,1);}\\
      {\tt r=r.*ndataset;}\\
      {\tt for j=1:ndataset}\\
      {\tt   i=floor(r(j));}\\
      {\tt   eval = c*1.0/(1.0+exp(-(a1*dataset\_in\_n(i)+a2*dataset\_in\_x(i)+b)))+d;}\\
      {\tt   error = 0.5*(eval-dataset\_y(i))$^2$;;}\\
      {\tt  derrorda1 = ( c*dataset\_in\_n(i)*exp(-a1*dataset\_in\_n(i)-a2*dataset\_in\_x(i)-b)*} \\
      {\tt \; \; \; (c / (exp(-a1*dataset\_in\_n(i)-a2*dataset\_in\_x(i)-b)+1)+d-dataset\_y(i))} \\
      {\tt \; \; \; ) / power((exp(-a1*dataset\_in\_n(i)-a2*dataset\_in\_x(i)-b)+1),2);} \\
      {\tt   a1=a1-eta*  derrorda1;}\\
      {\tt    derrorda2 = ( c*dataset\_in\_x(i)*exp(-a1*dataset\_in\_n(i)-a2*dataset\_in\_x(i)-b)*}\\
   {\tt  \; \; \;  (c / (exp(-a1*dataset\_in\_n(i)-a2*dataset\_in\_x(i)-b)+1)+d-dataset\_y(i))}\\
   {\tt  \; \; \;  ) / power((exp(-a1*dataset\_in\_n(i)-a2*dataset\_in\_x(i)-b)+1),2);}\\
  {\tt  a2=a2-eta*  derrorda2;}\\
  {\tt  derrordb = ( c*exp(-a1*dataset\_in\_n(i)-a2*dataset\_in\_x(i)-b)*}\\
  {\tt   \; \; \;   (c / (exp(-a1*dataset\_in\_n(i)-a2*dataset\_in\_x(i)-b)+1)+d-dataset\_y(i))}\\
  {\tt   \; \; \;  ) / power((exp(-a1*dataset\_in\_n(i)-a2*dataset\_in\_x(i)-b)+1),2);}\\
  {\tt  b=b-eta*  derrordb;}\\
  {\tt  derrordc = ( c / (exp(-a1*dataset\_in\_n(i)-a2*dataset\_in\_x(i)-b)+1)+d-dataset\_y(i)}\\
  {\tt   \; \; \;  ) / (exp(-a1*dataset\_in\_n(i)-a2*dataset\_in\_x(i)-b)+1);}\\
  {\tt  c=c-eta*  derrordc;}\\
  {\tt  derrordd = c / (exp(-a1*dataset\_in\_n(i)-a2*dataset\_in\_x(i)-b)+1)+d-dataset\_y(i);}\\
  {\tt  d=d-eta*  derrordd;}\\

      {\tt   \inred{\% evaluation of ANN approximation of sin(n*pi*x) for n=0.333}}\\
      {\tt   n=0.333;}\\
{\tt   x=0:0.01:0.5;}\\
{\tt   y=sin(n*pi.*x);}\\
{\tt   eval = c*1.0./(1.0+exp(-(a1*n+a2.*x+b)))+d;}\\
{\tt   plot(x,y,x,eval);}\\

\subsection{Verification}

In Figure \ref{t1a} we present the training over 250,000 samples, and in Figure \ref{t2a} we present the convergence of the training.

\begin{figure}[h]
\includegraphics[scale=0.4]{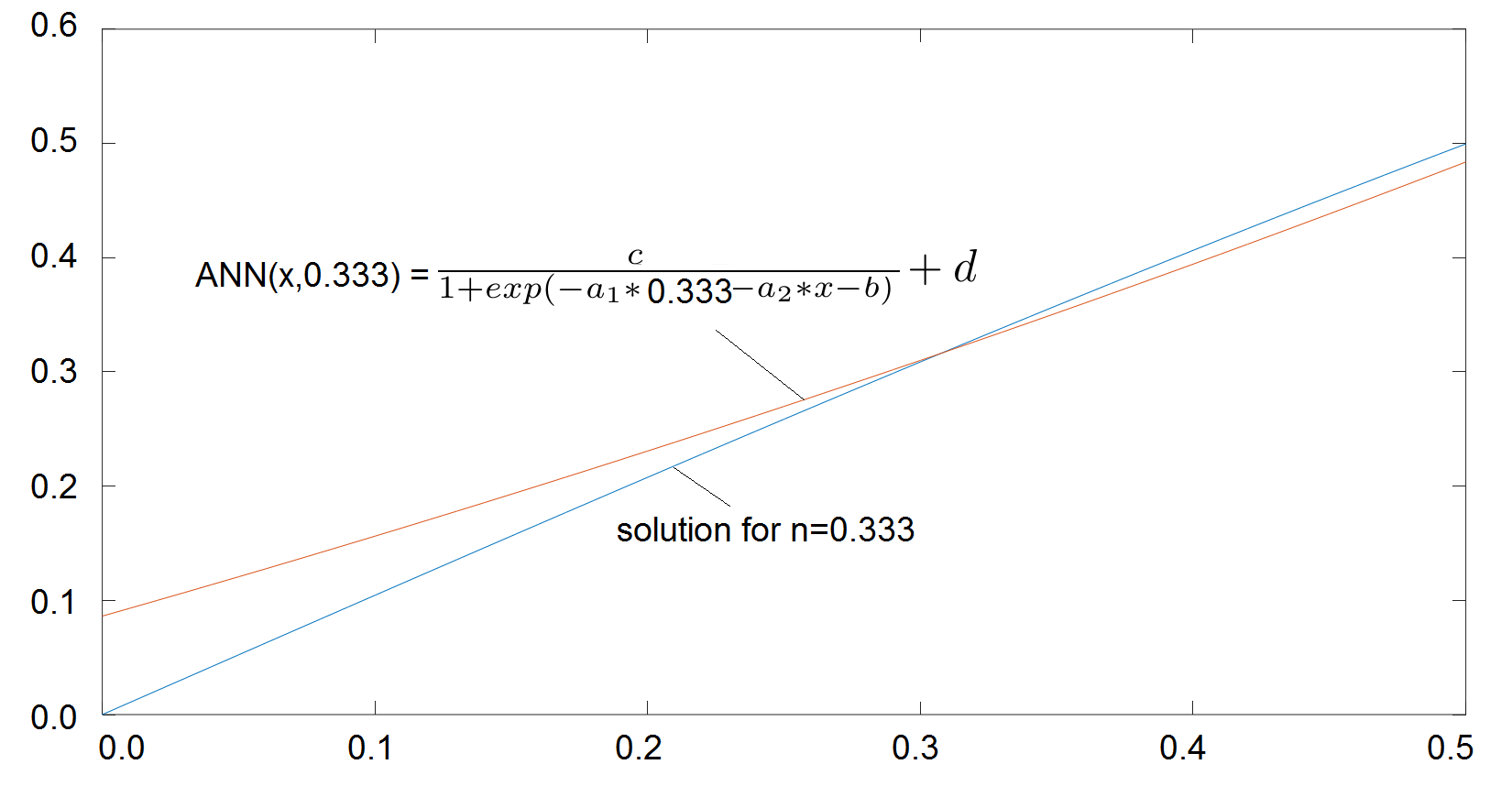}
\caption{Verification of the neural network approximation of solution for $n=0.333$  with
      $y(x)=ANN(n,x)= \frac{c}{1+exp(-a_1*n-a_2*x-b)}+d$
}\label{0333a}
\end{figure}

We select $n=0.333$ and we compute

\begin{eqnarray}
y(x)=ANN(n,x)= \frac{c}{1+exp(-a_1*n-a_2*x-b)}+d
\end{eqnarray}

we compare with $sin(0.333\pi x)$ in Figure \ref{0333a}.

\section{Appendix C: One-dimensional example of Physics Informed Neural Network }

\subsection{One dimensional heat-transfer problem}

We focus again on the heat-transfer problem
\begin{equation}
u''(x)+f(x)=0 \quad x \in (0,0.5)
\end{equation}
defined over $x \in (0,0.5)$,
with boundary conditions $u(0)=0$ and $u'(0.5)=g(0.5)$, with $g(x)=n \pi cos(n \pi x)$ and $f(x)=n^2 \pi^2 sin(n \pi x)$.

\subsection{Physics informed neural network}

We define the neural network

\begin{equation}
PINN(x)=u
\end{equation}
where
\begin{equation}
PINN(x)=c \sigma \left( a x +b \right)+d = \frac{c}{1+exp(-ax-b)}+d \label{ANN3}
\end{equation}

We compute the derivatives
\begin{equation}
PINN_x(x)=\frac{a*c*exp(-ax-b)}{\left(exp(-ax-b)+1\right)^2}
\end{equation}
and
\begin{equation}
PINN_{xx}(x)=c\left(\frac{2a^2 exp(-2ax-2b)}{\left(exp(-ax-b)+1\right)^3}-\frac{a^2 exp(-ax-b)}{\left(exp(-ax-b)+1\right)^2}\right)
\end{equation}

\subsection{Training}

The goal of the training is to find values of the weights $a,b,c,d$

We prepare a set of samples

\begin{itemize}
\item We randomly select $x \in (0,0.5)$
\item Input data $x$, output data $u=PINN(x)$
\end{itemize}

We define
\begin{equation}
F(x)=PINN_{xx}(x)+n^2 \pi^2 sin(n \pi x) 
\end{equation}

We define the error of approximation of PDE
\begin{eqnarray}
error1(x)=0.5*F(x)^2=0.5*\left(PINN_{xx}(x)+n^2 \pi^2 sin(n \pi x) \right)^2= \nonumber \\
0.5*\left(
c\left(\frac{2a^2 exp(-2ax-2b)}{\left(exp(-ax-b)+1\right)^3}-\frac{a^2 exp(-ax-b)}{\left(exp(-ax-b)+1\right)^2}\right)
+n^2 \pi^2 sin(n \pi x)
\right)^2
\end{eqnarray}
as well as the error of approximation of the boundary condition at $x=0$
\begin{eqnarray}
error2(0)=0.5*\left(PINN(0)-0\right)^2=0.5*\left(\frac{c}{1+exp(-b)}+d-0\right)^2
\end{eqnarray}
as well as the error of approximation of the boundary condition at $x=0.5$
\begin{eqnarray}
error3(0.5)=0.5*\left(PINN_x(0.5)-g(0.5)\right)^2= \nonumber \\ 
0.5*\left(\frac{a*c*exp(-a0.5-b)}{\left(exp(-a0.5-b)+1\right)^2}-n \pi cos(n \pi 0.5)\right)^2
\end{eqnarray}

\begin{enumerate}
\item Select $x$
\item Compute $u=PINN(x)=c \sigma \left( a x +b \right)+d = \frac{c}{1+exp(-ax-b)}+d$
\item Compute  $error1(x)$, $error2(0)$, $error3(0.5)$
\item Compute  $\frac{\partial error1(x)}{\partial a}, \frac{\partial error1(x)}{\partial b}, \frac{\partial error1(x)}{\partial c}, \frac{\partial error1(x)}{\partial d}$
\item Compute  $\frac{\partial error2(0)}{\partial a}, \frac{\partial error2(0)}{\partial b}, \frac{\partial error2(0)}{\partial c}, \frac{\partial error2(0)}{\partial d}$
\item Compute  $\frac{\partial error3(0.5)}{\partial a}, \frac{\partial error3(0.5)}{\partial b}, \frac{\partial error3(0.5)}{\partial c}, \frac{\partial error3(0.5)}{\partial d}$
\item Correct 
\begin{eqnarray}
a = a - \eta * \frac{\partial e(x)}{\partial a} \\
b = b - \eta * \frac{\partial e(x)}{\partial b} \\
c = c - \eta * \frac{\partial e(x)}{\partial c} \\
d = d - \eta * \frac{\partial e(x)}{\partial d} 
\end{eqnarray}
where $e(x)=error1(x)+error2(x)+error3(x)$
\end{enumerate}
for $\eta \in (0,1)$. 

We compute
\begin{eqnarray}
\frac{\partial PINN_{xx}(x)}{\partial a}=c\left(\frac{a^2x*exp(-ax-b)}{(exp(-ax-b)+1)^2}
-\frac{6a^2x*exp(-2ax-2b)}{(exp(-ax-b)+1)^3}\right. \nonumber \\
\left. +\frac{6a^2x*exp(-3ax-3b)}{(exp(-ax-b)+1)^4}
-\frac{2aexp(-ax-b)}{(exp(-ax-b)+1)^2}
+\frac{4a*exp(-2ax-2b)}{(exp(-ax-b)+1)^3}\right) 
\end{eqnarray}
\begin{eqnarray}
\frac{\partial PINN_{xx}(x)}{\partial b}=c\left(\frac{a^2exp(ax+b)\left(-4exp(ax+b)+exp(2ax+2b)+1\right)}{(exp(ax+b)+1)^4}\right)
\end{eqnarray}
\begin{equation}
\frac{\partial PINN_{xx}(x)}{\partial c}=\left(\frac{2a^2 exp(-2ax-2b)}{\left(exp(-ax-b)+1\right)^3}-\frac{a^2 exp(-ax-b)}{\left(exp(-ax-b)+1\right)^2}\right)
\end{equation}
\begin{equation}
\frac{\partial PINN_{xx}(x)}{\partial d}=0
\end{equation}
\begin{equation}
\frac{\partial PINN(0)}{\partial a}=0
\end{equation}
\begin{equation}
\frac{\partial PINN(0)}{\partial b}=\frac{exp(-b)c}{\left(exp(-b)+1\right)^2}
\end{equation}
\begin{equation}
\frac{\partial PINN(0)}{\partial c}=\frac{1}{\left(exp(-b)+1\right)}
\end{equation}
\begin{equation}
\frac{\partial PINN(0)}{\partial d}=1.0
\end{equation}
\begin{eqnarray}
\frac{\partial PINN_{x}(0.5)}{\partial a}=c*exp(b-a)\frac{\left((1-0.5a)*exp(2a+b)+(0.5a+1)exp(1.5a)\right)}{\left(exp(0.5a+b)+1\right)^3}
\end{eqnarray}
\begin{eqnarray}
\frac{\partial PINN_{x}(0.5)}{\partial b}=\frac{ac*exp(b-0.5a)\left(exp(a)-exp(1.5a+b)\right)}{\left(exp(0.5a+b)+1\right)^3}
\end{eqnarray}
\begin{eqnarray}
\frac{\partial PINN_{x}(0.5)}{\partial c}=\frac{a*exp(-0.5a-b)}{\left(exp(-0.5a-b)+1\right)^2}
\end{eqnarray}
\begin{eqnarray}
\frac{\partial PINN_{x}(0.5)}{\partial d}=0
\end{eqnarray}

Using the above formulas we have

\begin{eqnarray}
\frac{\partial error1(x)}{\partial a} = 
\left(c\left(\frac{2a^2 exp(-2ax-2b)}{\left(exp(-ax-b)+1\right)^3}-\frac{a^2 exp(-ax-b)}{\left(exp(-ax-b)+1\right)^2}\right)+n^2 \pi^2 sin(n \pi x) \right) \nonumber \\ c\left(\frac{a^2x*exp(-ax-b)}{(exp(-ax-b)+1)^2}
-\frac{6a^2x*exp(-2ax-2b)}{(exp(-ax-b)+1)^3}\right. \nonumber \\
\left. +\frac{6a^2x*exp(-3ax-3b)}{(exp(-ax-b)+1)^4}
-\frac{2aexp(-ax-b)}{(exp(-ax-b)+1)^2}
+\frac{4a*exp(-2ax-2b)}{(exp(-ax-b)+1)^3}\right) 
\end{eqnarray}
\begin{eqnarray}
\frac{\partial error1(x)}{\partial b} = 
\left(c\left(\frac{2a^2 exp(-2ax-2b)}{\left(exp(-ax-b)+1\right)^3}-\frac{a^2 exp(-ax-b)}{\left(exp(-ax-b)+1\right)^2}\right)+n^2 \pi^2 sin(n \pi x) \right) \nonumber \\ 
c\left(\frac{a^2exp(ax+b)\left(-4exp(ax+b)+exp(2ax+2b)+1\right)}{(exp(ax+b)+1)^4}\right)
\end{eqnarray}
\begin{eqnarray}
\frac{\partial error1(x)}{\partial c} = 
\left(c\left(\frac{2a^2 exp(-2ax-2b)}{\left(exp(-ax-b)+1\right)^3}-\frac{a^2 exp(-ax-b)}{\left(exp(-ax-b)+1\right)^2}\right)+n^2 \pi^2 sin(n \pi x) \right) \nonumber \\ 
\left(\frac{2a^2 exp(-2ax-2b)}{\left(exp(-ax-b)+1\right)^3}-\frac{a^2 exp(-ax-b)}{\left(exp(-ax-b)+1\right)^2}\right)
\end{eqnarray}
\begin{eqnarray}
\frac{\partial error1(x)}{\partial d} = 0
\end{eqnarray}

\begin{eqnarray}
\frac{\partial error2(x)}{\partial a} = 0
\end{eqnarray}
\begin{eqnarray}
\frac{\partial error2(x)}{\partial b} = \left(\frac{c}{1+exp(-b)}+d-0\right)
\frac{exp(-b)c}{\left(exp(-b)+1\right)^2}
\end{eqnarray}
\begin{eqnarray}
\frac{\partial error2(x)}{\partial c} = \left(\frac{c}{1+exp(-b)}+d-0\right) \nonumber \\
\left(\frac{2a^2 exp(-2ax-2b)}{\left(exp(-ax-b)+1\right)^3}-\frac{a^2 exp(-ax-b)}{\left(exp(-ax-b)+1\right)^2}\right)
\end{eqnarray}
\begin{eqnarray}
\frac{\partial error2(x)}{\partial d} = \left(\frac{c}{1+exp(-b)}+d-0\right)
\end{eqnarray}

\begin{eqnarray}
\frac{\partial error3(x)}{\partial a} = \left(\frac{a*c*exp(-a0.5-b)}{\left(exp(-a0.5-b)+1\right)^2}-n \pi cos(n \pi 0.5)\right)  \nonumber \\ 
c*exp(b-a)\frac{\left((1-0.5a)*exp(2a+b)+(0.5a+1)exp(1.5a)\right)}{\left(exp(0.5a+b)+1\right)^3}
\end{eqnarray}
\begin{eqnarray}
\frac{\partial error3(x)}{\partial b} = \left(\frac{a*c*exp(-a0.5-b)}{\left(exp(-a0.5-b)+1\right)^2}-n \pi cos(n \pi 0.5)\right) \nonumber \\ 
\frac{ac*exp(b-0.5a)\left(exp(a)-exp(1.5a+b)\right)}{\left(exp(0.5a+b)+1\right)^3}
\end{eqnarray}
\begin{eqnarray}
\frac{\partial error3(x)}{\partial c} = \left(\frac{a*c*exp(-a0.5-b)}{\left(exp(-a0.5-b)+1\right)^2}-n \pi cos(n \pi 0.5)\right) \nonumber \\ 
\frac{a*exp(-0.5a-b)}{\left(exp(-0.5a-b)+1\right)^2}
\end{eqnarray}
\begin{eqnarray}
\frac{\partial error3(x)}{\partial d} = 0
\end{eqnarray}

\begin{figure}[h]%
\includegraphics[scale=0.4]{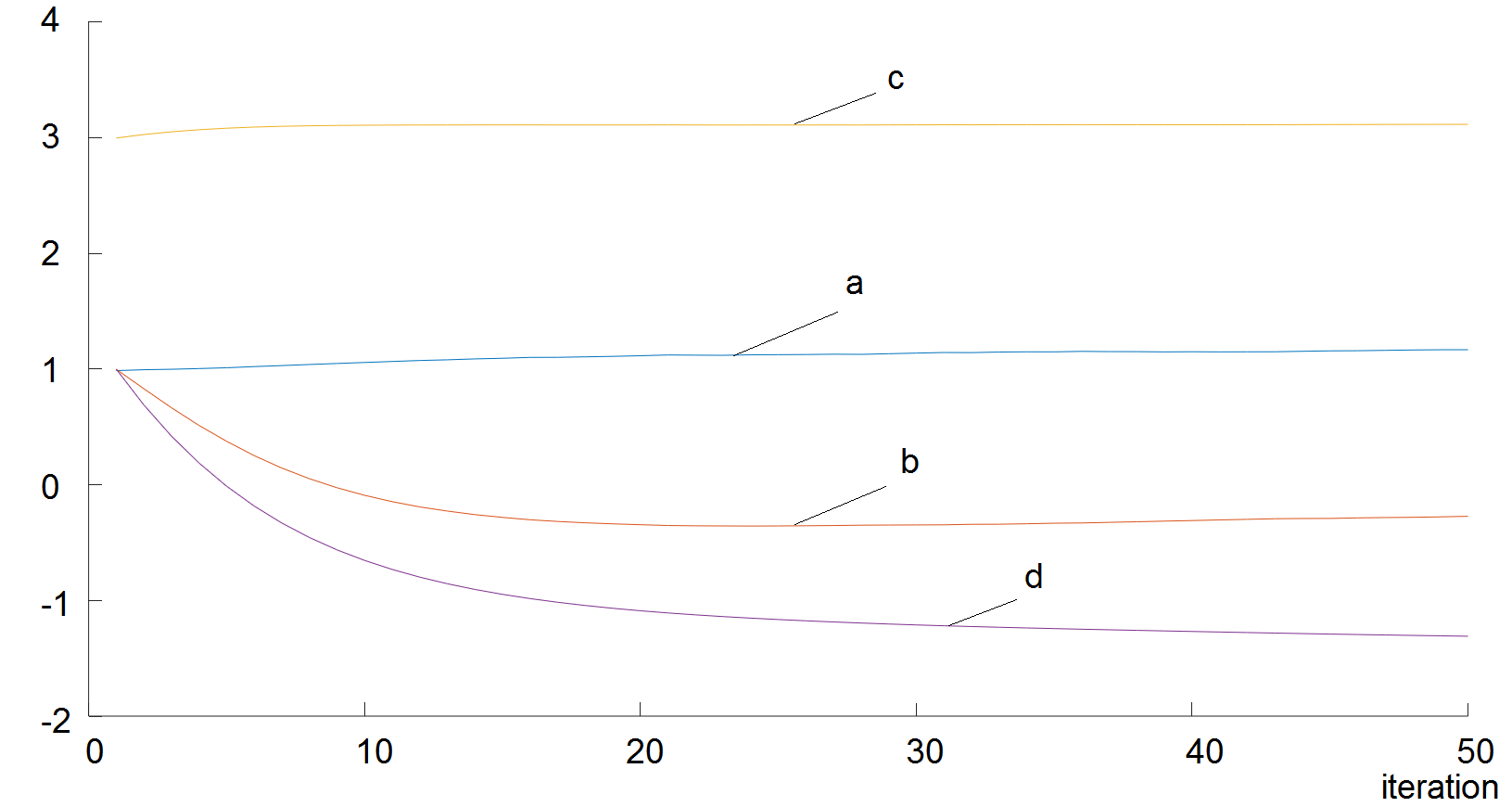}
\caption{Training of the simple PINN (\ref{ANN3})
 starting from $a=b=d=1.0$, and $c=3.0$, for $\eta=0.1$.} \label{t1b}
\end{figure}
\begin{figure}[h]
\includegraphics[scale=0.4]{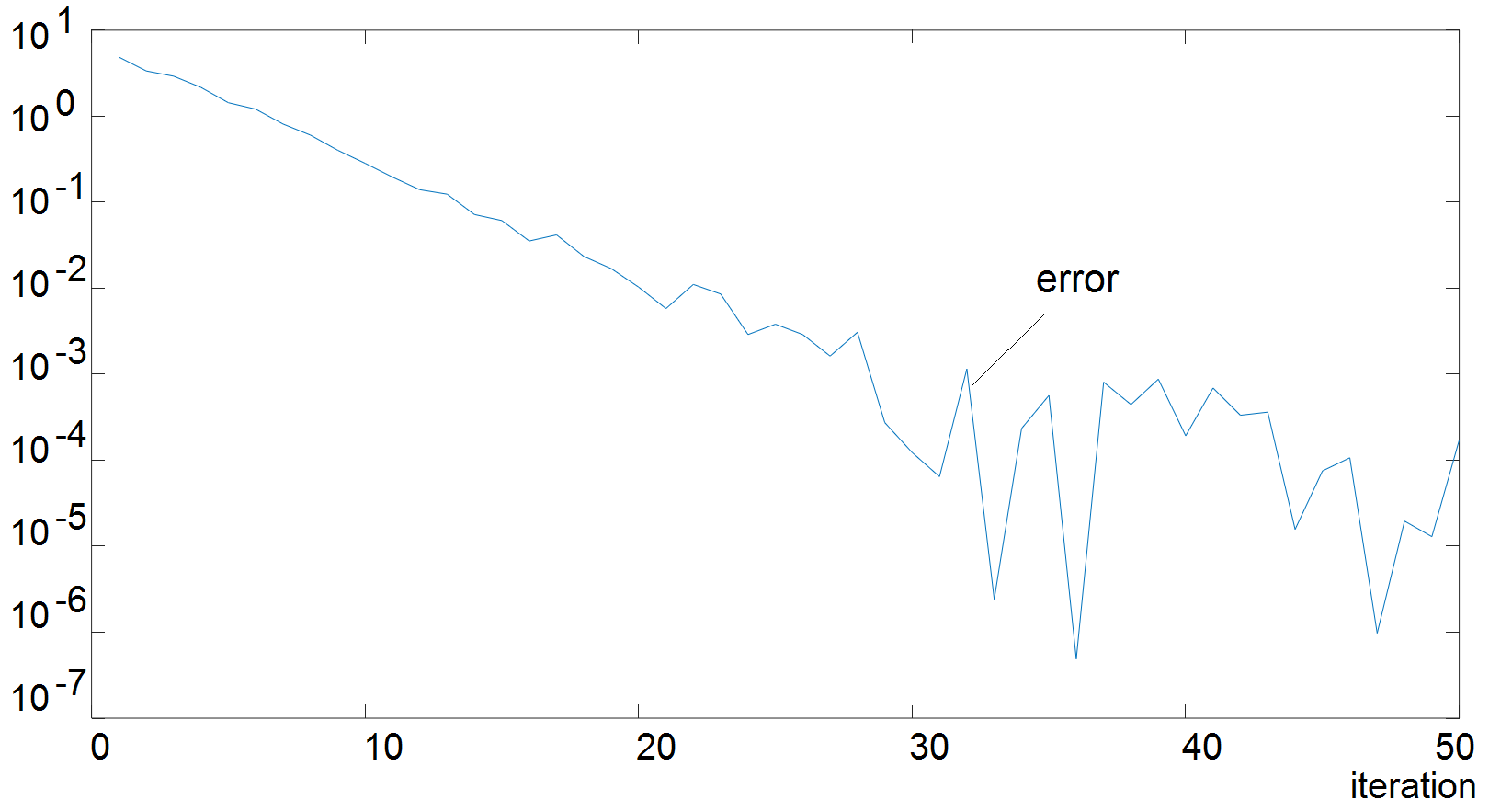}
\caption{Convergence of error for training of PINN} \label{t2b}
% starting from 1,1,10,-10}
\end{figure}

\subsection{MATLAB implementation}

      {\tt \inred{\% Creation of dataset}} \\
      {\tt i=1;} \\
      {\tt n=0.333;} \\
      {\tt for x=0.01:0.01:0.5} \\
      {\tt y=sin(n*pi*x);} \\
      {\tt dataset\_in\_x(i)=x; } \\
      {\tt dataset\_y(i)=y; } \\
      {\tt i=i+1;} \\
      {\tt endfor} \\
      {\tt ndataset=i-1;}\\

      {\tt \inred{\% Training}} \\
      {\tt a1=1.0; a2=1.0; b=1.0; c=3.0; d=1.0;}\\
      {\tt eta=0.1;}\\
      {\tt r = 0 + (1-0).*rand(ndataset,1);}\\
      {\tt r=r.*ndataset;}\\
      {\tt for j=1:ndataset}\\
      {\tt   i=floor(r(j));}\\
      {\tt   eval = c*1.0/(1.0+exp(-(a*dataset\_in\_x(i)+b)))+d;}\\
      {\tt   error = 0.5*(eval-dataset\_y(i))$^2$;}\\
      {\tt   \inred{\% Training of the PDE}}\\
      {\tt   x = dataset\_in\_x(i); }\\
	{\tt    n = dataset\_in\_n(i);}\\
	{\tt      Fx = c*(2*a*a*exp(-2*a*x-2*b) / power((exp(-a*x-b)+1),3) }\\
      {\tt     - a*a*exp(-a*x-b) / power((exp(-a*x-b)+1),2))+n*n*pi*pi*sin(n*x);}\\
	{\tt      derror1da = Fx*c *( (a*a*x*exp(-a*x-b))/power((exp(-a*x-b)+1),2)  }\\
	{\tt      -(6*a*a*x*exp(-2*a*x-2*b))/power((exp(-a*x-b)+1),3) }\\
	{\tt      +(6*a*a*x*exp(-3*a*x-3*b))/power((exp(-a*x-b)+1),4) }\\ 
	{\tt      -(2*a*exp(-a*x-b))/power((exp(-a*x-b)+1),2) }\\
	{\tt      +(4*a*exp(-2*a*x-2*b))/power((exp(-a*x-b)+1),3) );}\\
	{\tt      a=a-eta*  derror1da;}\\
	{\tt      derror1db =  Fx*c*( (a*a*exp(a*x+b))*(-4*exp(a*x+b) }\\
	{\tt      +exp(2*a*x+2*b)+1)/power((exp(a*x+b)+1),4));}\\
	{\tt      b=b-eta*  derror1db;}\\
	{\tt      derror1dc = Fx*( (2*a*a*exp(-2*a*x-2*b))/power((exp(-a*x-b)+1),3)}\\
	{\tt      -(a*a*exp(-a*x-b))/power((exp(-a*x-b)+1),2)); }\\
	{\tt      c=c-eta*  derror1dc; }\\
	{\tt      derror1dd = 0; }\\
	{\tt      d=d-eta*  derror1dd; }\\
      {\tt   \inred{\% Training of the boundary condition at x=0}}\\
	{\tt        x=0; }\\
	{\tt        derror2da = 0;}\\
	{\tt        a=a-eta*  derror2da;}\\
	{\tt        derror2db = (c / (1+exp(-b))+d)* (exp(-b)*c)/power((exp(-b)+1),2);}\\
	{\tt        b=b-eta*  derror2db;}\\
	{\tt        derror2dc = (c/(1+exp(-b))+d)*( (2*a*a*exp(-2*a*x-2*b))/power(exp(-a*x-b)+1,3)}\\
     {\tt  - (a*a*exp(-a*x-b))/power(exp(-a*x-b)+1,2) );}\\
	{\tt        c=c-eta*  derror2dc;}\\
	{\tt        derror2dd = c/(1+exp(-b))+d;  }\\
	{\tt        d=d-eta*  derror2dd;}\\
      {\tt   \inred{\% Training of the boundary condition at x=0.5}}\\
	{\tt        x=0.5;}\\
	{\tt        derror3da = (a*c*exp(-a*0.5-b)/power((exp(-a*0.5-b)+1),2)}\\
     {\tt -n*pi*cos(n*pi*0.5)) *c*exp(b-a)*((1-0.5*a)*exp(2*a+b)}\\
     {\tt +(0.5*a+1)*exp(1.5*a))/power((exp(0.5*a+b)+1),3);}\\
	{\tt        a=a-eta*  derror3da;}\\
	{\tt        derror3db = (a*c*exp(-a*0.5-b)/power((exp(-a*0.5-b)+1),2)}\\
      {\tt -n*pi*cos(n*pi*0.5))* a*c*exp(b-0.5*a)*(exp(a)+exp(1.5*a+b))/power((exp(0.5*a+b)+1),3);}\\
	{\tt        b=b-eta*  derror3db;}\\
	{\tt        derror3dc = (a*c*exp(-a*0.5-b)/power((exp(-a*0.5-b)+1),2)}\\
     {\tt -n*pi*cos(n*pi*0.5))* a*exp(-b-0.5*a)/power((exp(-0.5*a-b)+1),2);}\\
	{\tt        c=c-eta*  derror3dc;}\\
	{\tt        derror3dd = 0;}\\
	{\tt        d=d-eta*  derror3dd; }\\
{\tt endfor }\\

      {\tt   \inred{\% evaluation of PINN approximation of sin(0.333*pi*x)}}\\
      {\tt   n=0.333;}\\
{\tt   x=0:0.01:0.5;}\\
{\tt   y=sin(n*pi.*x);}\\
{\tt   eval = c*1.0./(1.0+exp(-(a.*x+b)))+d;}\\
{\tt   plot(x,y,x,eval);}\\

\subsection{Verification}

In Figure \ref{t1b} we present the training over 50 samples, and in Figure \ref{t2b} we present the convergence of the training.

\begin{figure}[h]
\includegraphics[scale=0.4]{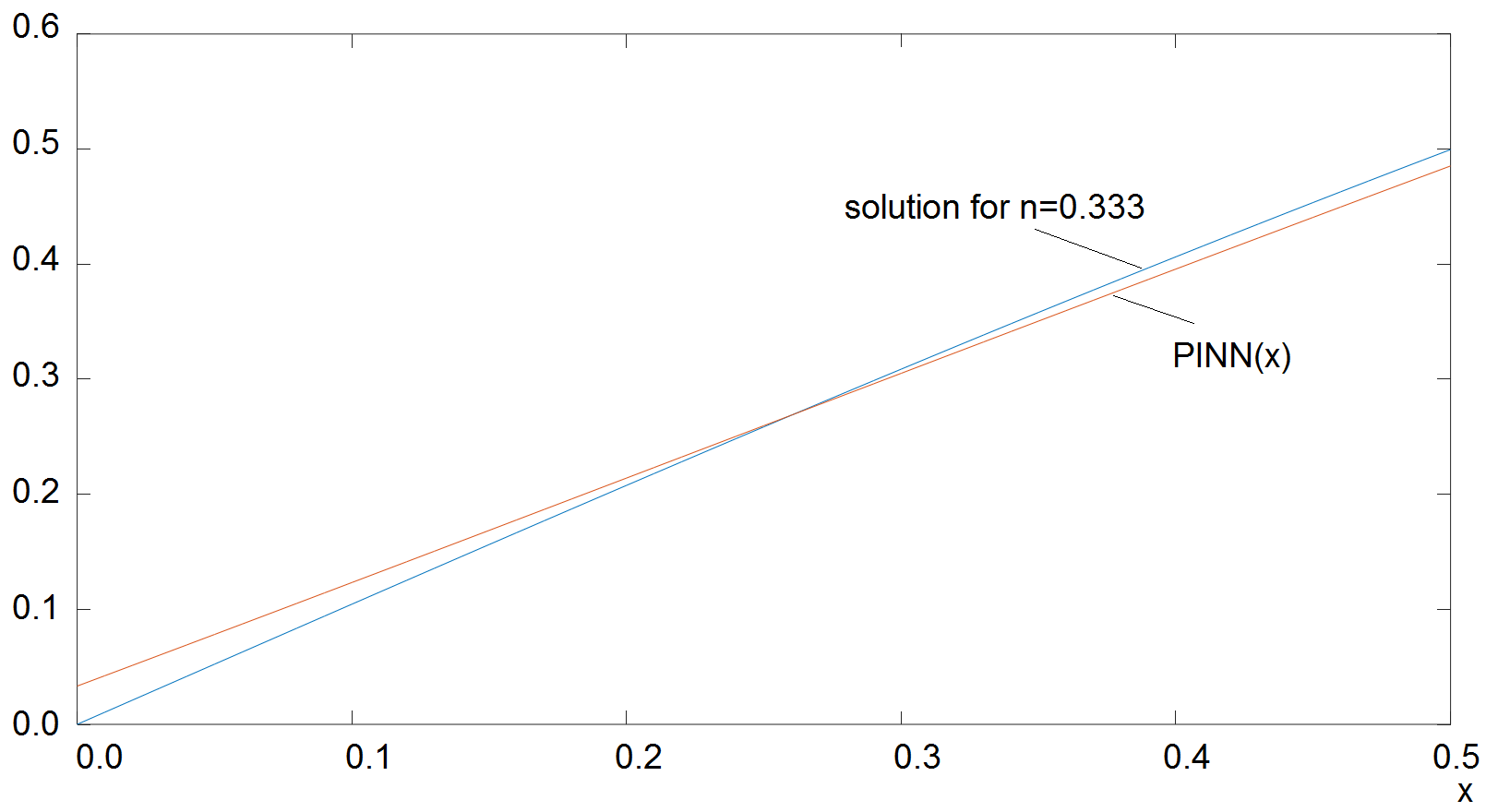}
\caption{Verification of the PINN trained for $n=0.333$  with
      $y(x)=PINN(n,x)= \frac{c}{1+exp(a*x-b)}+d$
}\label{0333b}
\end{figure}

The PINN has been trained for $n=0.333$ so we compute

\begin{eqnarray}
y(x)=ANN(n,x)= \frac{c}{1+exp(a*x-b)}+d
\end{eqnarray}

we compare with $sin(0.333\pi x)$ in Figure \ref{0333b}.

\bibliographystyle{named}
\bibliography{biblio}

\begin{thebibliography}{}

\bibitem[\protect\citeauthoryear{Brevis \bgroup \em et al.\egroup
  }{2021}]{BREVIS2021186}
Ignacio Brevis, Ignacio Muga, and Kristoffer~G. {van der Zee}.
\newblock A machine-learning minimal-residual (ml-mres) framework for
  goal-oriented finite element discretizations.
\newblock {\em Computers \& Mathematics with Applications}, 95:186--199, 2021.
\newblock Recent Advances in Least-Squares and Discontinuous Petrov–Galerkin
  Finite Element Methods.

\bibitem[\protect\citeauthoryear{Chen \bgroup \em et al.\egroup }{2018}]{NODE}
Ricky T.~Q. Chen, Yulia Rubanova, Jesse Bettencourt, and David Duvenaud.
\newblock Neural ordinary differential equations.
\newblock {\em Advances in Neural Information Processing Systems 31}, 2018.

\bibitem[\protect\citeauthoryear{Collier \bgroup \em et al.\egroup
  }{2012}]{CostOfContinuity}
Nathan Collier, David Pardo, Lisandro Dalcin, Maciej Paszynski, and V.M. Calo.
\newblock The cost of continuity: A study of the performance of isogeometric
  finite elements using direct solvers.
\newblock {\em Computer Methods in Applied Mechanics and Engineering},
  213-216:353--361, 2012.

\bibitem[\protect\citeauthoryear{Dalcin \bgroup \em et al.\egroup
  }{2016}]{PetIGA}
L.~Dalcin, N.~Collier, P.~Vignal, A.M.A. Côrtes, and V.M. Calo.
\newblock Petiga: A framework for high-performance isogeometric analysis.
\newblock {\em Computer Methods in Applied Mechanics and Engineering}, 308:151
  -- 181, 2016.

\bibitem[\protect\citeauthoryear{Haghighat \bgroup \em et al.\egroup
  }{2020}]{PINNIGA}
Ehsan Haghighat, Maziar Raissi, Adrian Mourec, Hector Gomez, and Ruben Juanes.
\newblock A deep learning framework for solution and discovery in solid
  mechanics.
\newblock {\em arXiv:2003.02751}, 2020.

\bibitem[\protect\citeauthoryear{Hornik \bgroup \em et al.\egroup
  }{1989}]{UniversalApproximators}
Kurt Hornik, Maxwell Stinchcombe, and Halbert White.
\newblock Multilayer feedforward networks are universal approximators.
\newblock {\em Neural Networks}, 1989.

\bibitem[\protect\citeauthoryear{HP}{2021}]{HP}
Deep learning driven self-adaptive hp finite element method.
\newblock {\em Lecture Notes in Computer Science}, (12742):114--121, 2021.

\bibitem[\protect\citeauthoryear{Hughes \bgroup \em et al.\egroup
  }{2005}]{IsogeometricAnalysisProposal}
Thomas~J.R. Hughes, J.A. Cottrell, and Yuri Bazilevs.
\newblock Isogeometric analysis: Cad, finite elements, nurbs, exact geometry
  and mesh refinement.
\newblock {\em Computer Methods in Applied Mechanics and Engineering}, 2005.

\bibitem[\protect\citeauthoryear{Lagaris \bgroup \em et al.\egroup }{1998}]{bc}
I.E. Lagaris, A.~Likas, and D.I. Fotiadis.
\newblock Artificial neural networks for solving ordinary and partial
  differential equations.
\newblock {\em IEEE Transactions on Neural Networks}, 9(5):987--1000, 1998.

\bibitem[\protect\citeauthoryear{Michoski \bgroup \em et al.\egroup
  }{2020}]{DiffEqDNN}
Craig Michoski, Miloš Milosavljević, Todd Oliver, and David~R. Hatch.
\newblock Solving differential equations using deep neural networks.
\newblock {\em Neurocomputing}, 399:193--212, 2020.

\bibitem[\protect\citeauthoryear{Misyris \bgroup \em et al.\egroup
  }{2019}]{PINNForPowerSystems}
George~S. Misyris, Andreas Venzke, and Spyros Chatzivasileiadis.
\newblock Physics-informed neural networks for power systems.
\newblock preprint at
  https://www.researchgate.net/publication/337184659\_Physics-Informed\_Neural\_Networks\_for\_Power\_Systems,
  2019.

\bibitem[\protect\citeauthoryear{Nguyen \bgroup \em et al.\egroup
  }{2015}]{IGAOvwCompImpl}
Vinh~Phu Nguyen, Cosmin Anitescu, Stéphane~P.A. Bordas, and Timon Rabczuk.
\newblock Isogeometric analysis: An overview and computer implementation
  aspects.
\newblock {\em Mathematics and Computers in Simulation}, 117:89 -- 116, 2015.

\bibitem[\protect\citeauthoryear{Pang \bgroup \em et al.\egroup
  }{2020}]{nPINNs}
Guofei Pang, Marta D'Elia, Michael Parks, and George Karniadakis.
\newblock npinns: Nonlocal physics-informed neural networks for a parametrized
  nonlocal universal laplacian operator. algorithms and applications.
\newblock {\em Journal of Computational Physics}, 422, 2020.

\bibitem[\protect\citeauthoryear{Paszy\'nski}{2020}]{IzogeometrycznaMES}
Maciej Paszy\'nski.
\newblock {\em Classical and isogeometric finite element method}.
\newblock AGH University of Science and Technology in Krakow, 2020.
\newblock
  https://epodreczniki.open.agh.edu.pl/handbook/1088/module/1173/reader.

\bibitem[\protect\citeauthoryear{Raissi \bgroup \em et al.\egroup
  }{2017}]{RaissiPhysicsIDL}
M.~Raissi, P.~Perdikaris, and G.~Karniadakis.
\newblock Physics informed deep learning (part i): Data-driven solutions of
  nonlinear partial differential equations.
\newblock {\em ArXiv}, abs/1711.10561, 2017.

\bibitem[\protect\citeauthoryear{Raissi \bgroup \em et al.\egroup
  }{2019}]{PDESolving}
Maziar Raissi, Paris Perdikaris, and George~E Karniadakis.
\newblock Physics-informed neural networks: A deep learning framework for
  solving forward and inverse problems involving nonlinear partial differential
  equations.
\newblock {\em Journal of Computational Physics}, 378:686--707, 2019.

\bibitem[\protect\citeauthoryear{Wang \bgroup \em et al.\egroup
  }{2020}]{WhenPinnsFail}
Sifan Wang, Xinling Yu, and Paris Perdikaris.
\newblock When and why pinns fail to train: A neural tangent kernel
  perspective.
\newblock preprint at \url{https://arxiv.org/abs/2007.14527}, 07 2020.

\end{thebibliography}

\end{document}